\documentclass[12pt]{article}
\usepackage[utf8]{inputenc}
\usepackage[english]{babel}
\usepackage{graphicx}
\usepackage{vmumath}
\usepackage{amsmath}
\usepackage{amsthm}
\usepackage{avo_book}
\usepackage{tikz}
\usepackage{caption}
\usepackage{subcaption}
\usepackage{booktabs}
\usepackage[a4paper,left=15mm,right=15mm,top=30mm,bottom=20mm]{geometry}
\noaftermath
\begin{document}

\bibliographystyle{unsrt}

\title{Enumeration of Chord Diagrams without Loops and Parallel Chords}

\author{
Evgeniy Krasko \qquad  Alexander Omelchenko\\
\small St. Petersburg Academic University\\
\small 8/3 Khlopina Street, St. Petersburg, 194021, Russia\\
\small\tt \{krasko.evgeniy, avo.travel\}@gmail.com
}


\begin{abstract}
We enumerate chord diagrams without loops and without both loops and parallel chords. We show that the former ones describe Hamiltonian paths in $n$-dimensional octahedrons. The latter ones are also known as shapes. For labelled diagrams we obtain generating functions, for unlabelled ones we derive recurrence relations.

\end{abstract}

\maketitle

\section{Introduction}

\begin{figure}[ht]
\centering
	\centering
    	\includegraphics[scale=0.8]{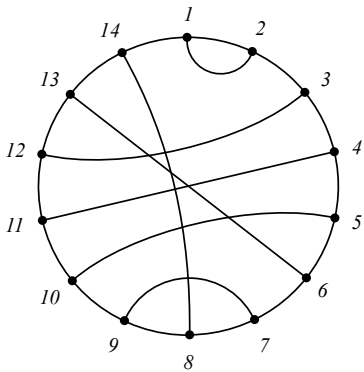}
	\caption{A chord diagram}
\label{fig:chord}
\end{figure}

A chord diagram consists of $2n$ points on a circle labelled with the numbers $1,2,\ldots,2n$ in a circular order, joined pairwise by chords (figure \ref{fig:chord}). Many different combinatorial objects are equivalent to chord diagrams, for instance one-face maps \cite{Walsh_Lehman} or polygons with sides identified pairwise \cite{HarerZagier}. Chord diagrams appear in theory of finite type invariants of knots and links \cite{Bar-Natan}, the representation theory of Lie algebras \cite{Campoamor}, the geometry of moduli spaces of flat connections on surfaces \cite{And_Mat_Reshet}, in mapping class groups \cite{Andersen2} and in other fields of pure mathematics. They also find applications in RNA analysis \cite{RNA}, \cite{Penner}, classification of textile structures \cite{textile_1}, \cite{textile_2} and analysis of data structures \cite{Flajolet}. 

In practice, different classes of chord diagrams appear. In this paper we study two of them, which we will call loopless diagrams and simple diagrams.  A chord is said to be a loop if it connects two neighboring points (chord $\{1,2\}$ on figure \ref{fig:chord}). Two chords are called parallel if they connect two pairs of neighboring points and don't intersect. For example, chords $\{3,12\}$ and $\{4, 11\}$ on figure \ref{fig:chord} are parallel, as well as chords $\{4,11\}$ and $\{5, 10\}$, but chords $\{3,12\}$ and $\{5, 10\}$ are not: it would be convenient for us that being parallel is not a transitive relation. 

A loopless chord diagram is just a diagram without loops. It turns out that such diagrams correspond to Hamiltonian cycles in an $n$-dimensional octahedron. Consequently, by enumerating loopless chord diagrams we enumerate these cycles as well. 

 By a simple diagram we mean a diagram that doesn't have neither loops nor pairs of parallel chords. Simple diagrams correspond to so-called shapes, e.g. one-face maps without vertices of degree $1$ or $2$  \cite{Chapuy_Schaeffer}. They are very important to the analysis of arbitrary one-face maps, as they describe their basic structure in the following sense: any one-face map can be obtained from a unique shape by subdividing its edges and attaching trees to it  \cite{Chapuy_Schaeffer}.

Depending on the notion of isomorphism used, two diagrams are said to be isomorphic if one could be obtained from the other either by a rotation or by a combination of rotations and reflections of the circle. Isomorphism classes of labelled chord diagrams are said to be unlabelled ones. Such diagrams find wide application in generating knots and links on surfaces \cite{k-tangle-projections} and for calculating knot invariants \cite{Kauffman_textile}. 

Labelled loopless chord diagrams were enumerated in \cite{Kalashnikov} using a combinatorial approach. There exist formulas for the numbers of unlabelled chord diagrams without restrictions \cite{Sawada}, \cite{Khruzin}. Unlabelled loopless chord diagrams were implicitly counted up to a small number of chords with the help of a computer in \cite{Singmaster} as Hamiltonian paths in an $n$-dimensional octahedron.

In this paper we focus on enumerating loopless and simple chord diagrams. For both classes of labelled diagrams we obtain explicit generating functions. We also derive a multivariate generating function that classifies chord diagrams according to the numbers of loops and pairs of parallel chords. As an intermediate result, we provide generating functions for the corresponding classes of linear diagrams, which are not circles but segments with points identified pairwise. 

For unlabelled objects we derive systems of recurrences that can be used to efficiently compute the numbers of loopless and simple chord diagrams. In both cases we give the answers for two kinds of symmetries: only rotations as well as rotations and reflections. Finally, we provide tables listing the numbers obtained from those generating functions and recurrences for chord diagrams with up to 20 chords.

\section{Labelled loopless chord diagrams}

\begin{figure}[ht]
\centering
	\begin{subfigure}[b]{0.26\textwidth}
	\centering
    		\includegraphics[scale=0.8]{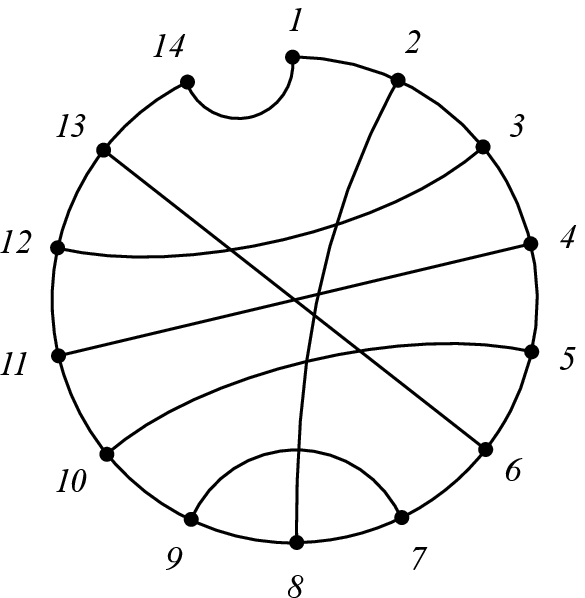}
		\caption{Cutting a chord diagram}
	\end{subfigure}
	\begin{subfigure}[b]{0.71\textwidth}
	\centering
    		\includegraphics[scale=0.8]{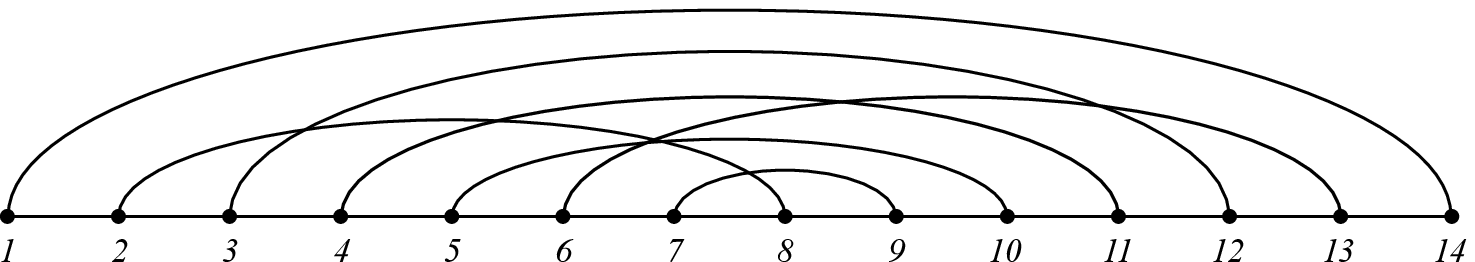}
		\caption{The same diagram drawn on a line}
	\end{subfigure}	
	\caption{A linear diagram}
\label{fig:linear}
\end{figure}

The number of all chord diagrams having $n$ chords is easy to calculate. Indeed, to obtain a chord diagram we could connect the point $1$ to any of $2n-1$ other points, the next free point with $2n-3$ remaining ones and so on. Thus the number of chord diagrams is equal to $(2n-1)!!$. The exponential generating function $b(t)$ for these numbers has the following form: 
$$
b(t)=\dfrac{1}{\sqrt{1-2t}}=\sum\limits_{n=0}^{\infty}(2n-1)!!\dfrac{t^n}{n!}.
$$
For loopless and simple diagrams it would be convenient to study so-called linear diagrams together with chord ones. Imagine we cut the circle of some chord diagram between points $1$ and $2n$. As a result we obtain a linear diagram, which is different from the original chord diagram in that points $1$ and $2n$ are no longer considered to be neighboring (see figure \ref{fig:linear}). For loopless diagrams we have a simple connection between the numbers of linear ($a_n$) and chord ($b_n$) diagrams:
\begin{equation}
\label{eq:dd_b_nk}
b_n=a_n-a_{n-1},\qquad n\geq 2;\qquad b_1=0.
\end{equation}
Indeed, if we glue a loopless linear diagram with $n\geq 2$ chords into a chord diagram, that is, start thinking of points $1$ and $2n$ as of neighbors, we get a diagram with a loop if and only if these points were initially connected with a chord. The number of linear diagrams with this chord is obviously equal to the number $a_{n-1}$ of linear loopless diagrams with $n-1$ chords, as the operation of removing the chord $(1,2n)$ together with its endpoints establishes a bijection between the corresponding sets. 

\begin{figure}[ht]
\centering
	\begin{subfigure}[b]{0.3\textwidth}
	\centering
    		\includegraphics[scale=0.8]{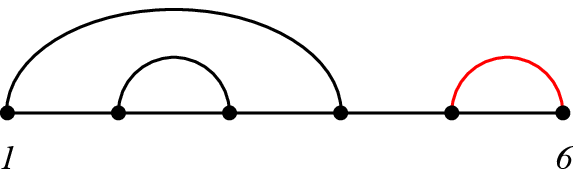}
		\caption{}
	\end{subfigure}
	\begin{subfigure}[b]{0.3\textwidth}
	\centering
    		\includegraphics[scale=0.8]{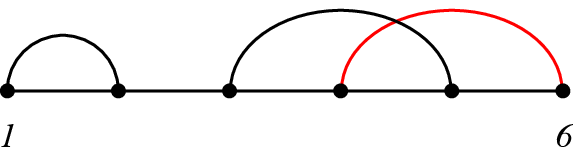}
		\caption{}
	\end{subfigure}
	\begin{subfigure}[b]{0.3\textwidth}
	\centering
    		\includegraphics[scale=0.8]{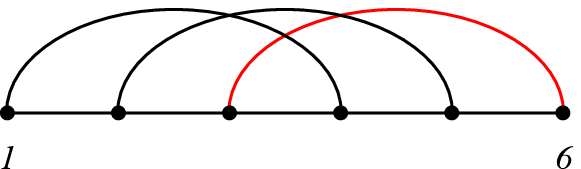}
		\caption{}
	\end{subfigure}	
	\caption{Counting linear diagrams}
\label{fig:loops_linear}
\end{figure}

So we reduced the problem of enumerating loopless chord diagrams to the enumeration of loopless linear ones. Instead of enumerating them directly, it would be convenient for us to derive a recurrence for the numbers $a_{n,k}$ of linear diagrams with $n$ chords, $k$ of which are loops. Let us show that
\begin{equation}
\label{eq:dd_a_nk}
a_{n+1,k}=a_{n,k-1}+(2n-k)a_{n,k}+(k+1)a_{n,k+1},\qquad a_{n,k}=0\quad \text{if $k>n$ or $k<0$,}\qquad a_{0,0}=1.
\end{equation} 
Consider a chord $(i,2n+2)$ in a diagram with $n+1$ chords and $k$ loops. There could be three different cases. If $i=2n+1$ removing this chord yields a linear diagram with $k-1$ loops (figure~\ref{fig:loops_linear}(a)). This gives us the summand $a_{n,k-1}$. Diagrams corresponding to the summand $(k+1)a_{n,k-1}$ are built from diagrams with $k+1$ loops by inserting a point in the middle of one of these loops and connecting it with another new point added to the right of the rightmost point of the initial diagram (figure \ref{fig:loops_linear}(b)). The multiplier $(k+1)$ is explained by the fact that the loop could be chosen in $k+1$ ways. 

It remains to explain the summand $(2n-k)a_{n,k}$ in (\ref{eq:dd_a_nk}). The corresponding diagrams are built in such a way: take the diagram with $k$ loops and $n$ chords and add a new point so that it does not break any loops (figure \ref{fig:loops_linear}(c)). Then connect it with another new point added to the right of the diagram. The multiplier $2n-k$  comes from the fact that the initial diagram had $2n-1$ intervals between consecutive points and one more position to the left of it. However, $k$ of these positions are forbidden for the first new point, because they correspond to loops.

For solving the recurrence (\ref{eq:dd_a_nk}) we introduce the polynomials
$$
P_n(z)=a_{n,0}+a_{n,1}z+a_{n,2}z^2+\ldots+a_{n,n}z^n,
$$
multiply (\ref{eq:dd_a_nk}) by $z^k$ and sum over $k$ from $0$ to $n+1$:
$$
\sum\limits_{k=0}^{n+1}a_{n+1,k}z^k=\sum\limits_{k=0}^{n+1}a_{n,k-1}z^k+2n\sum\limits_{k=0}^{n+1}a_{n,k}z^k-
\sum\limits_{k=0}^{n+1}ka_{n,k}z^k+ \sum\limits_{k=0}^{n+1}(k+1)a_{n,k+1}z^k.
$$
Taking into account the boundary conditions $a_{n,k}=0$ if $k>n$ or $k<0$, we get the following equation for the function $P_n(z)$:
\begin{equation}
\label{eq:bred_0}
P_{n+1}(z)=zP_{n}(z)+2nP_{n}(z)-zP'_{n}(z)+P'_{n}(z)=zP_{n}(z)+2nP_{n}(z)+(1-z)P'_{n}(z).
\end{equation}
The initial condition $a_{0,0}=1$ in terms of $P_n(z)$ could be rewritten as $P_0(z)=1$. To solve the infinite system of equations (\ref{eq:bred_0}) it is convenient to define the following two-variable generating function
$$
w(z,t)=\sum\limits_{n=0}^{+\infty}P_n(z)\dfrac{t^{n}}{n!}.
$$
Multiplying (\ref{eq:bred_0}) by $t^n/n!$ and summing over $n$ from $0$ to $+\infty$, we get
\begin{equation}
\label{eq:bred_1}
\sum\limits_{n=0}^{+\infty}P_{n+1}(z)\dfrac{t^{n}}{n!}=z\sum\limits_{n=0}^{+\infty}P_{n}(z)\dfrac{t^{n}}{n!}+
2\sum\limits_{n=1}^{+\infty}P_n(z)\dfrac{t^n}{(n-1)!}+(1-z)\sum\limits_{n=0}^{+\infty}P'_n(z)\dfrac{t^n}{n!},
\end{equation}
which could be rewritten in terms of $w(z,t)$ as
$$
\pd{w}{t}=zw(z,t)+2t\pd{w}{t}+(1-z)\pd{w}{z}.
$$
The condition $P_0(z)=1$ takes the form $w(z,0)=1$. The solution for this Cauchy problem is
$$
w(z,t)=\dfrac{e^{(-1+\sqrt{1-2t})(1-z)}}{\sqrt{1-2t}}.
$$
Now the generating function $\phi(t)$ for loopless linear diagrams is easily obtained by setting $z=0$ in $w(z,t)$:
$$
\phi(t)=w(0,t)=\dfrac{e^{-1+\sqrt{1-2t}}}{\sqrt{1-2t}}=1+0\cdot t+1\cdot\dfrac{t^2}{2!}+5\cdot\dfrac{t^3}{3!}+36\cdot\dfrac{t^4}{4!}+
329\cdot\dfrac{t^5}{5!}+\ldots
$$
This function allows us to derive some useful recurrences for the numbers $a_n$. As an example, we express the derivative
$$
\phi'(t)=\dfrac{e^{-1+\sqrt{1-2t}}}{(1-2t)^{3/2}}\left[1-\sqrt{1-2t}\right]
$$
through the function $\phi(t)$:
$$
\phi'(t)=\phi(t)\left[\dfrac{1-\sqrt{1-2t}}{1-2t}\right]\qquad\Longleftrightarrow\qquad
(1-2t)\phi'(t)=\phi(t)\left[1-\sqrt{1-2t}\right].
$$
This yields the following recurrence:
$$
a_{n+1}=2na_n+\sum\limits_{k=1}^n\BCf{n}{k}(2k-3)!!a_{n-k}.
$$
Using the second derivative
$$
\phi''(t)=\dfrac{e^{-1+\sqrt{1-2t}}}{\sqrt{1-2t}}\left[\dfrac{1-2t+3[1-\sqrt{1-2t}]}{(1-2t)^2}\right]=
\dfrac{\phi(t)+3\phi'(t)}{1-2t},
$$
we could get from the equation
$$
\phi''(t)[1-2t]=\phi(t)+3\phi'(t)
$$
the following second-order recurrence:
\begin{equation}
\label{eq:three}
a_{n+2}=(2n+3)a_{n+1}+a_n\qquad\Longleftrightarrow\qquad a_{n+1}=(2n+1)a_n+a_{n-1};\qquad a_0=1,\quad a_1=0.
\end{equation}
According to \cite{Kalashnikov}, it was first guessed by Jean Betramas of LABRI, Bordeaux, from the numerical computations and then combinatorially proved by Michiel Hazewinkel and V.~V.~Kalashnikov \cite{Kalashnikov}.

Now we turn back to the equality (\ref{eq:dd_b_nk}) which connects the numbers $b_n$ and $a_n$. Let $\psi(t)=\sum_{n=1}^{+\infty}b_n t^n/n!$ be the exponential generating function for the numbers $b_n$. Multiplying (\ref{eq:dd_b_nk}) by $t^n/n!$ and summing the result over $n$, we get 
\begin{equation}
\label{eq:bn_an_sum}
\sum\limits_{n=2}^{\infty}b_n\dfrac{t^n}{n!}= \sum\limits_{n=2}^{\infty}a_n\dfrac{t^n}{n!}-
\sum\limits_{n=2}^{\infty}a_{n-1}\dfrac{t^n}{n!}\qquad \Longleftrightarrow\qquad
\psi(t)=\phi(t)-a_0-a_1t-\chi(t),
\end{equation}
where $\chi(t)$ is the generation function for the numbers $a_{n-1}$. This function could be expressed through the integral of $\phi(t)$:
$$
\chi(t)=\sum\limits_{n=2}^{\infty}a_{n-1}\dfrac{t^n}{n!}=\sum\limits_{n=1}^{\infty}a_{n}\dfrac{t^{n+1}}{(n+1)!}=
\int\phi(t)\,dt-a_0t=1-t-e^{-1+\sqrt{1-2t}}.
$$
It follows that the generation function for labelled loopless chord diagrams is equal to
$$
\psi(t)=\phi(t)-2+t+e^{-1+\sqrt{1-2t}}=e^{-1+\sqrt{1-2t}}\left(1+\dfrac{1}{\sqrt{1-2t}}\right)-2+t=
$$
$$
=0\cdot t+1\cdot\dfrac{t^2}{2!}+4\cdot\dfrac{t^3}{3!}+31\cdot\dfrac{t^4}{4!}+
293\cdot\dfrac{t^5}{5!}+3326\cdot\dfrac{t^6}{6!}+\ldots
$$
(sequence $A003436$ in oeis.org).

\section{Unlabelled loopless chord diagrams}

To count the numbers $\tilde{b}_n$ of unlabelled loopless chord diagrams with $n$ chords we will use the Burnside lemma
\begin{equation}
\label{eq:Burnside_lemma}
\tilde{b}_n=\dfrac{1}{|G|}\sum\limits_{g\in G}|\Fix(g)|.
\end{equation}
Here $|\Fix(g)|$ is the number of labelled diagrams fixed by the action of the element $g$ of some group $G$ that defines the isomorphism relation between diagrams. In our case $G$ will be either the cyclic group $C_{2n}$ of diagram's rotations or the dihedral group $D_{2n}$ of rotations and reflections. 

We start with the simpler case of the cyclic group $C_{2n}$. Consider the action of the group $C_{2n}$ on the set of loopless chord diagrams with $2n$ points and $n$ chords. Let $d$ be a divisor of $2n$, $\phi(d)$ be the Euler function of it. There are $\phi(d)$ elements of order $d$ in $C_{2n}$. Any such element fixes the same number $f(2n,d)$ of diagrams. These diagrams we will call $d$-symmetric. So, (\ref{eq:Burnside_lemma}) could be rewritten as
\begin{equation}
\label{eq:Burnside_lemma_Cn}
\tilde{b}_n=\dfrac{1}{2n}\sum\limits_{d\,|\, 2n}\phi(d)\,f(2n,d).
\end{equation} 

\begin{figure}[ht]
\centering
	\begin{subfigure}[b]{0.4\textwidth}
	\centering
    		\includegraphics[scale=0.8]{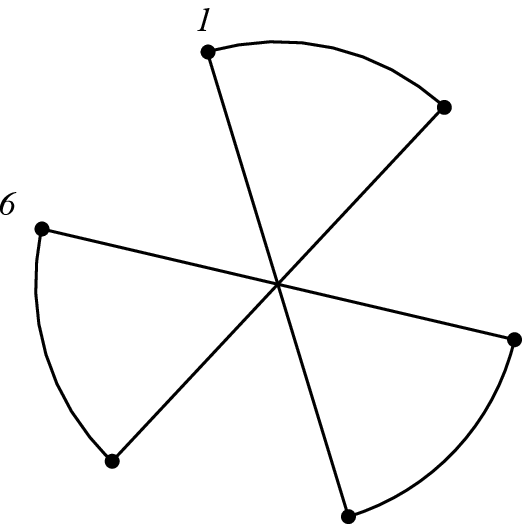}
		\caption{}
	\end{subfigure}
	\begin{subfigure}[b]{0.4\textwidth}
	\centering
    		\includegraphics[scale=0.8]{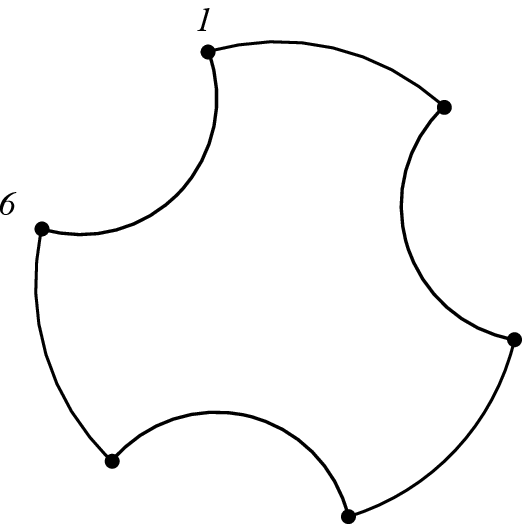}
		\caption{}
	\end{subfigure}	
	\caption{3-symmetric linear diagrams}
\label{fig:d_symmetric}
\end{figure}

Instead of enumerating $d$-symmetric chord diagrams directly we start with counting so-called $d$-symmetric linear diagrams (figure \ref{fig:d_symmetric}(a)). Any such a diagram with $2n$ points is obtained by cutting the circle of a $d$-symmetric chord diagram into $d$ sectors, each of them having $m:=2n/d$ points, between points $m$ and $m+1$, $2m$ and $2m+1$, $\ldots$, $2n$ and $1$. By cutting between points $i$ and $i+1$ we mean again that these points are no longer considered to be neighbors. Note that $1$-symmetric linear diagrams are just linear diagrams considered in the previous section. If we cut a $d$-symmetric loopless chord diagram, we get a loopless $d$-symmetric linear diagram. The converse is not true: if points $m$ and $m+1$, $2m$ and $2m+1$, $\ldots$, $2n$ and $1$ are connected by chords in a loopless linear $d$-symmetric diagram, then gluing this diagram back into a chord diagram results in $d$ loops in it (figure \ref{fig:d_symmetric}(b)). 

\begin{figure}[ht]
\centering
	\begin{subfigure}[b]{0.24\textwidth}
	\centering
    		\includegraphics[scale=0.8]{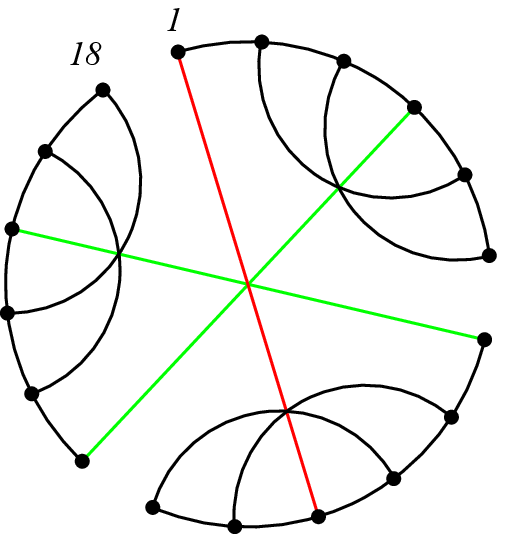}
		\caption{}
	\end{subfigure}	
	\begin{subfigure}[b]{0.24\textwidth}
	\centering
    		\includegraphics[scale=0.8]{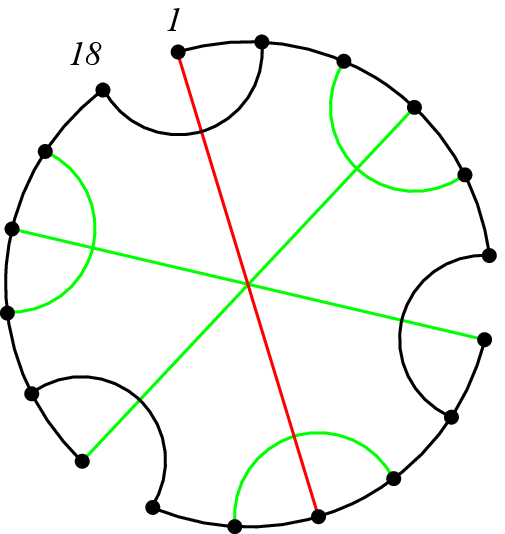}
		\caption{}
	\end{subfigure}
	\begin{subfigure}[b]{0.24\textwidth}
	\centering
    		\includegraphics[scale=0.8]{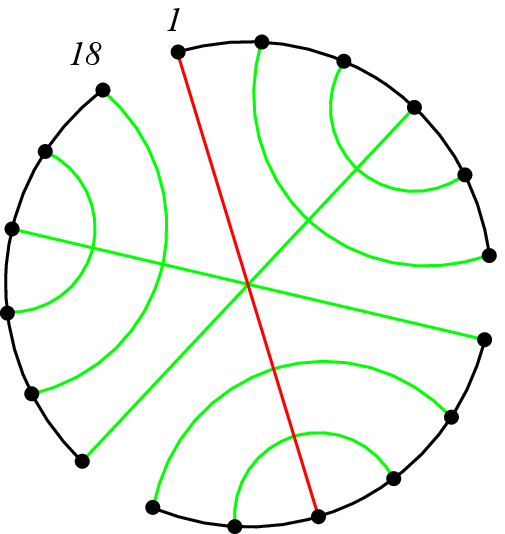}
		\caption{}
	\end{subfigure}	
	\begin{subfigure}[b]{0.24\textwidth}
	\centering
    		\includegraphics[scale=0.8]{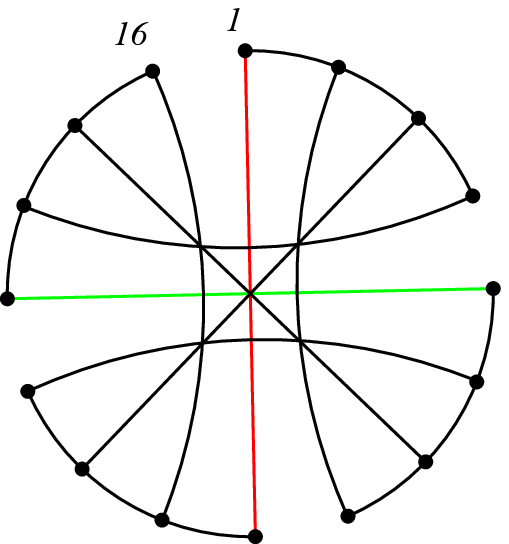}
		\caption{}
	\end{subfigure}
	\caption{$d$-symmetric loopless linear diagrams}
\label{fig:loopless_chord_d_symmetric_odd}
\end{figure}

Denote by $a_m^{(d)}$ the number of $d$-symmetric linear diagrams having $m\cdot d$ points. Our next goal is to derive a recurrence relation for these numbers. Consider first the case of odd $d$. Let us remove the chord $(1,i)$ together with its endpoints as well all the other chords on its orbit (under a rotation with a period $d$) in a $d$-symmetric linear diagram. In the simple case the result of this operation is a loopless $d$-symmetric linear diagram with $(m-2)\cdot d$ points. This is illustrated on figure \ref{fig:loopless_chord_d_symmetric_odd}(a). Hereafter, the red chord is the one that we remove first, and the green ones are those that we are obliged to remove subsequently to preserve some diagram's property like symmetry or looplessness. 

In a more complex case we get $d$ loops in the resulting diagram (figure \ref{fig:loopless_chord_d_symmetric_odd}(b)). If that happens, we remove these loops too. However we could get $d$ new loops again (figure \ref{fig:loopless_chord_d_symmetric_odd}(c)) and so on. As a result we obtain the following recurrence for the numbers $a_m^{(d)}$:
\begin{equation}
\label{eq:a_n_odd_0}
a_m^{(d)}=[d(m-1)-1]\cdot a_{m-2}^{(d)}+\sum_{i=1}^{m/2-1}d\cdot (m-1-2i)\cdot a_{m-2-2i}^{(d)}.
\end{equation}
The multiplier $d(m-1)-1$ in (\ref{eq:a_n_odd_0}) is explained by the fact that the second endpoint $i$ of the chord $(1,i)$ can be any point in $\{1,\ldots,2n\}$ except for the points $1,2,m,2m,3m,\ldots,(d-1)m$. The multipliers $d\cdot [m-1-2i]$ could be explained by similar arguments.

Rewriting the equality (\ref{eq:a_n_odd_0}) for $a_{m-2}^{(d)}$ and subtracting it from (\ref{eq:a_n_odd_0}) we obtain 
\begin{equation}
\label{eq:a_n_odd}
a_m^{(d)}=d(m-1)a_{m-2}^{(d)}+a_{m-4}^{(d)}.
\end{equation}
Observe that by setting $d=1$ we prove (\ref{eq:three}) combinatorially in slightly different terms.  

Note that in the case of odd $d$ we always delete a multiplier of $d$ chords at a time. In the case of even $d$ it could happen that we should delete only $d/2$ chords. This corresponds to the case of a chord $(1,n+1)$ connecting the opposite points of a diagram (figure \ref{fig:loopless_chord_d_symmetric_odd}(d)). As a consequence one more summand $a_{m-1}^{(d)}$ is added to the recurrence relation:
\begin{equation}
\label{eq:a_n_even_0}
a_m^{(d)}=a_{m-1}^{(d)}+[d(m-1)-1]\cdot a_{m-2}^{(d)}+\sum_{i=1}^{m/2-1}d\cdot (m-1-2i)\cdot a_{m-2-2i}^{(d)}.
\end{equation}
Rewriting (\ref{eq:a_n_even_0}) for $a_{m-2}^{(d)}$ and subtracting from (\ref{eq:a_n_even_0}) as before, we get
\begin{equation}
\label{eq:a_n_even}
a_m^{(d)}=a_{m-1}^{(d)}+d(m-1)a_{m-2}^{(d)}-a_{m-3}^{(d)}+a_{m-4}^{(d)}.
\end{equation}
The initial conditions for (\ref{eq:a_n_odd}) and (\ref{eq:a_n_even}) are the following:
$$
a_m^{(d)}=0,\quad m < 0; \qquad a_0^{(d)} = 1; \qquad a_{1}^{(2k)}=1; \qquad  a_{1}^{(2k+1)}=0; \qquad a_{2}^{(2k)}=d; \qquad  a_{2}^{(2k+1)}=d-1.
$$

\begin{figure}[ht]
	\centering
    	\includegraphics[scale=0.8]{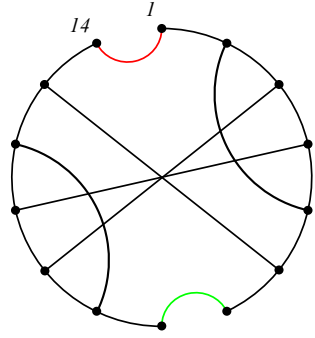}
	\caption{Special case of a $d$-symmetric loopless linear diagram}
\label{fig:loopless_chord}
\end{figure}

Now we get the expression for the number $f(m\cdot d,d)$ of $d$-symmetric loopless chord diagrams with $m\cdot d> 2$ points:
$$
f(m\cdot d,d)=a_m^{(d)}-a_{m-2}^{(d)}.
$$
Indeed, any loopless $d$-symmetric linear diagram with $m\cdot d>2$ points corresponds to a loopless $d$-symmetric chord diagram if and only if it has no chord $(m,m+1)$. The number of diagrams having such a chord is clearly equal to $a_{m-2}^{(d)}$ (figure \ref{fig:loopless_chord}). If $m\cdot d=2$ the number $b_m^{(d)}$ is equal to zero. 

Substituting the last expression into (\ref{eq:Burnside_lemma_Cn}), we obtain the numbers $\tilde{b}_n$ of loopless chord diagrams not isomorphic under rotations: $0,1,2,7,36,300,3218,\ldots$ (see Table \ref{table:loopless}).

\begin{figure}[ht]
\centering
	\begin{subfigure}[b]{0.3\textwidth}
	\centering
    		\includegraphics[scale=0.8]{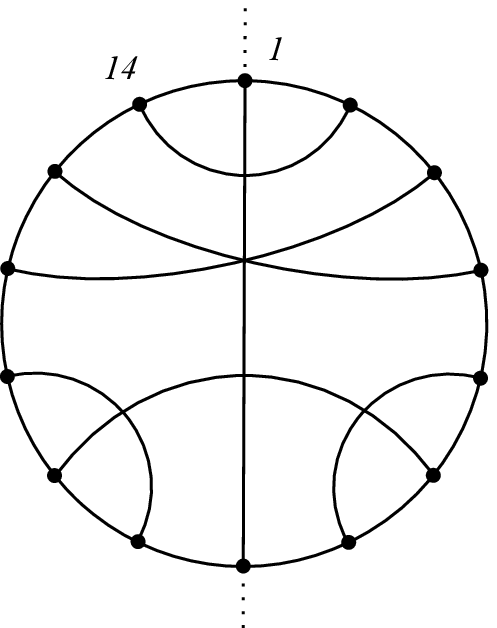}
		\caption{}
	\end{subfigure}
	\begin{subfigure}[b]{0.3\textwidth}
	\centering
    		\includegraphics[scale=0.8]{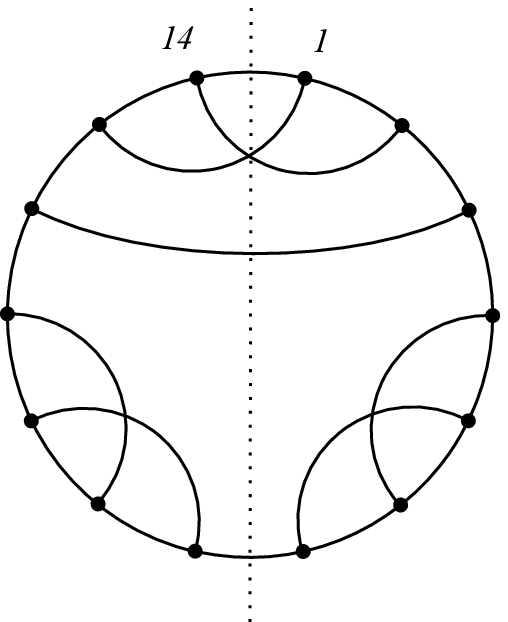}
		\caption{}
	\end{subfigure}
	\begin{subfigure}[b]{0.3\textwidth}
	\centering
    		\includegraphics[scale=0.8]{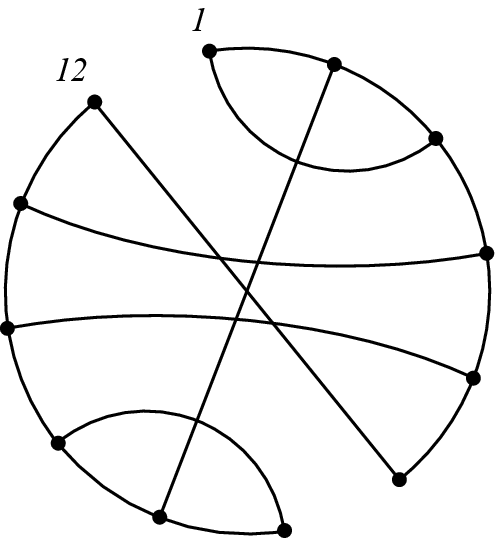}
	\caption{}
	\end{subfigure}	
	\begin{subfigure}[b]{0.3\textwidth}
	\centering
    		\includegraphics[scale=0.8]{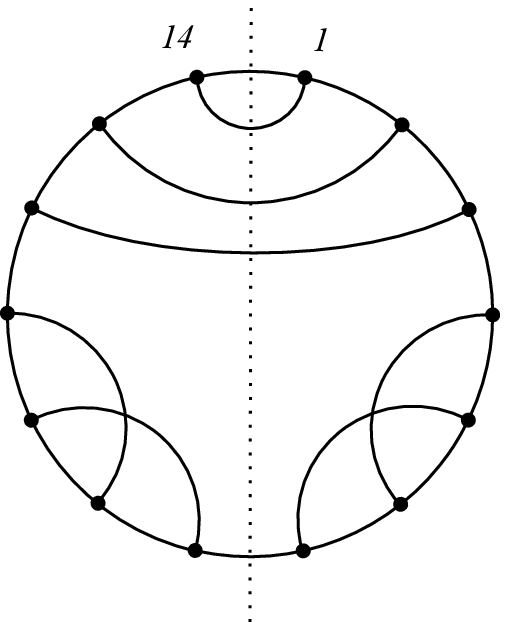}
	\caption{}
	\end{subfigure}	
	\begin{subfigure}[b]{0.3\textwidth}
	\centering
    		\includegraphics[scale=0.8]{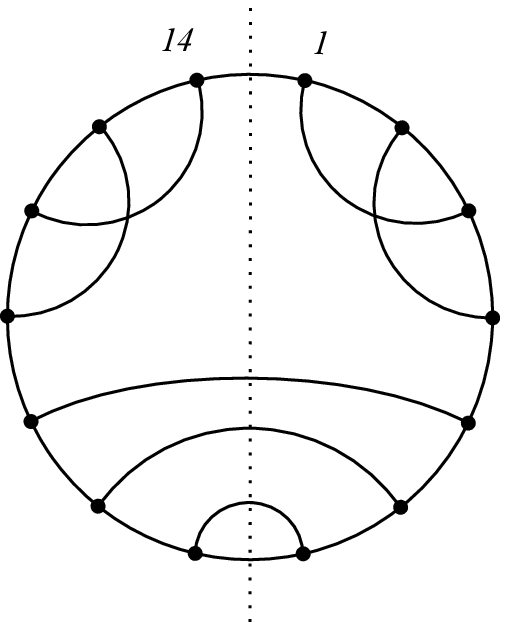}
	\caption{}
	\end{subfigure}
	\begin{subfigure}[b]{0.3\textwidth}
	\centering
    		\includegraphics[scale=0.8]{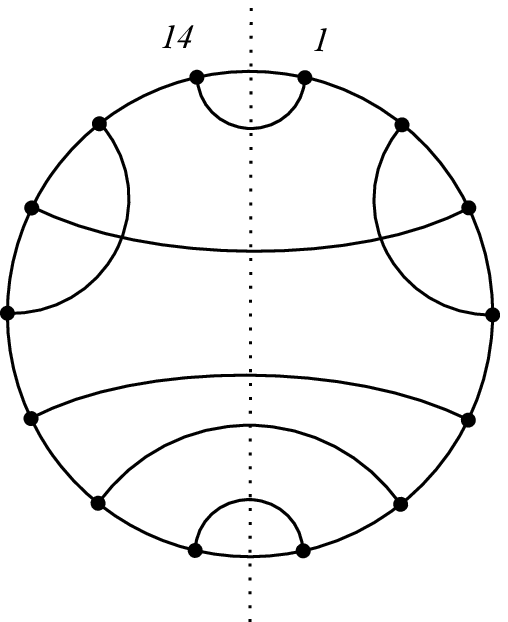}
	\caption{}
	\end{subfigure}	
	\caption{Diagrams with reflectional symmetry}
\label{fig:loopless_chord_refl}
\end{figure}

We now proceed with enumerating non-isomorphic chord diagrams under the action of the dihedral group $D_{2n}$. The Burnside lemma could be rewritten for this case as
\begin{equation}
\label{eq:Burnside_lemma_Dn}
\tilde{c}_n=\dfrac{1}{4n}\left[\sum\limits_{d\,|\, 2n}\phi(d)\,f(2n,d)+n\cdot K(n)+n\cdot H(n)\right],
\end{equation} 
where $K(n)$ denotes the number of chord diagrams with $n$ chords symmetric under the reflection about the axis passing through two opposite points of a diagram (figure \ref{fig:loopless_chord_refl}(a)), $H(n)$ is the number of diagrams that are symmetric under the reflection about the axis passing through the midpoints of the arcs connecting neighboring points of a diagram (figure \ref{fig:loopless_chord_refl}(b)). 

We claim that the numbers $K(n)$ are equal to $a_{n-1}^{(2)}$. Indeed, any diagram of the corresponding kind could be obtained by taking a $2$-symmetric loopless linear diagram, flipping its right half, and adding one more chord lying of the axis of symmetry (figure \ref{fig:loopless_chord_refl}(a) can be obtained from figure \ref{fig:loopless_chord_refl}(c) in such a way). We also claim that
$$
H(n)=a_{n}^{(2)}-2a_{n-1}^{(2)}+a_{n-2}^{(2)},\qquad n \geq 2;  \qquad\qquad H(2)=0.
$$
To prove this equality let us take any $2$-symmetric loopless linear diagram with $n$ chords and flip its right half as before. Gluing it into a chord diagram results in not more than $2$ loops. The number of loopless diagrams can be obtained with the help of inclusion-exclusion principle: take the number $a_n^{(2)}$ of all diagrams, subtract twice the number $a_{n-1}^{(2)}$ of diagrams with at least $1$ loop (figure \ref{fig:loopless_chord_refl}(d) and (e)) and add the number $a_{n-2}^{(2)}$ of diagrams with exactly $2$ loops (figure \ref{fig:loopless_chord_refl}(f)). 

Finally we get the following formula for the numbers of loopless chord diagrams:
$$
\tilde{c}_n=\dfrac{\tilde{b}_n}{2}+\dfrac{a_{n}^{(2)}-a_{n-1}^{(2)}+a_{n-2}^{(2)}}{4},\qquad  n \geq 2;\qquad\qquad 
\tilde{c}_1=0,\quad \tilde{c}_2=1.
$$
The first terms of the corresponding sequence $0,1,2,7,29,1788,\ldots$ (see Table \ref{table:loopless}) were first obtained by D. Singmaster in \cite{Singmaster} as a result of numerical computation of non-isomorphic Hamiltonian paths in an $n$-dimensional octahedron (the sequence $A003437$ in oeis.org). Let us show that there indeed exists a bijection between Hamiltonian paths in octahedrons and loopless chord diagrams.

\begin{figure}[ht]
\centering
	\begin{subfigure}[b]{0.3\textwidth}
	\centering
    		\includegraphics[scale=0.8]{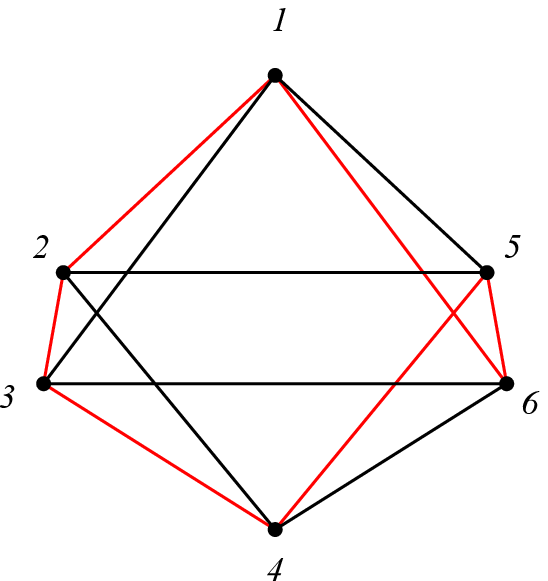}
		\caption{}
	\end{subfigure}	
	\begin{subfigure}[b]{0.3\textwidth}
	\centering
    		\includegraphics[scale=0.8]{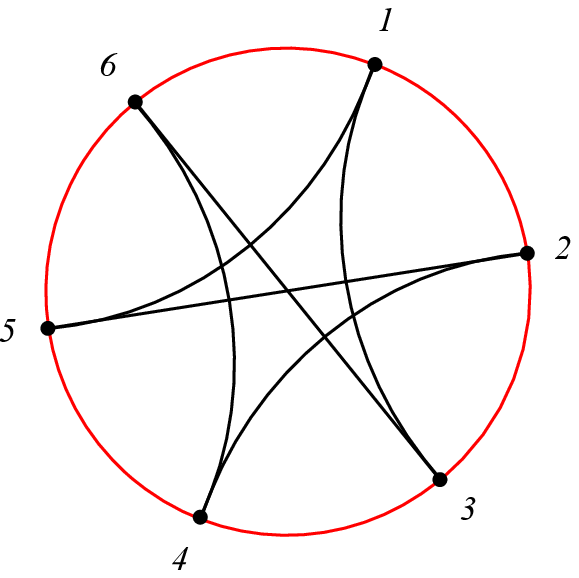}
		\caption{}
	\end{subfigure}
	\begin{subfigure}[b]{0.3\textwidth}
	\centering
    		\includegraphics[scale=0.8]{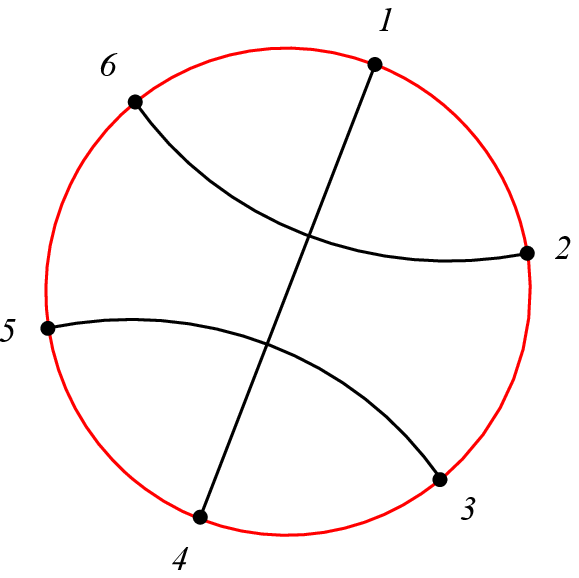}
		\caption{}
	\end{subfigure}
	\caption{Correspondence between Hamiltonian cycles in octahedrons and chord diagrams}
\label{fig:octahedron}
\end{figure}

An $n$-dimensional octahedron is a graph having $n$ pairs of vertices; every vertex is connected to all the rest except for the one from its pair. Take an $n$-dimensional octahedron with a distinguished Hamiltonian path (figure \ref{fig:octahedron}(a)) and draw it in such a way that this path forms a circle in a plane (figure \ref{fig:octahedron}(b)). Now remove all of its edges not belonging to the Hamiltonian path and add chords between those vertices that were not connected by an edge (figure \ref{fig:octahedron}(c)). The resulting object is a chord diagram which is necessary loopless: loops would correspond to parallel edges in an octahedron. Clearly, this transformation is invertible. 

Labels on figure \ref{fig:octahedron} are added just to visualize the correspondence between vertices. Any unlabelled octahedron with a distinguished unlabelled and undirected Hamiltonian path corresponds to an unlabelled loopless chord diagram considered with respect to both rotational and reflectional symmetries in exactly the same way.

\section{Labelled and unlabelled simple chord diagrams}

\begin{figure}[ht]
\centering
	\begin{subfigure}[b]{0.32\textwidth}
	\centering
    		\includegraphics[scale=0.8]{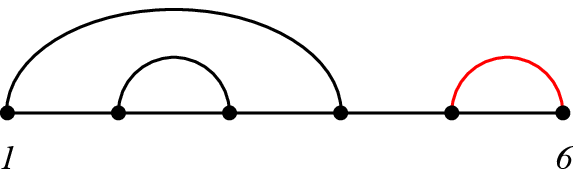}
		\caption{}
	\end{subfigure}	
	\begin{subfigure}[b]{0.32\textwidth}
	\centering
    		\includegraphics[scale=0.8]{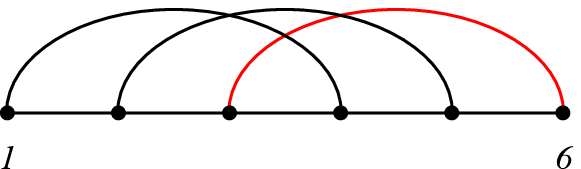}
		\caption{}
	\end{subfigure}
	\begin{subfigure}[b]{0.32\textwidth}
	\centering
    		\includegraphics[scale=0.8]{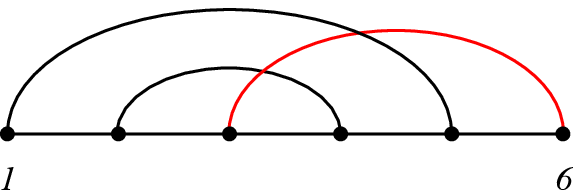}
		\caption{}
	\end{subfigure}
	\begin{subfigure}[b]{0.45\textwidth}
	\centering
    		\includegraphics[scale=0.8]{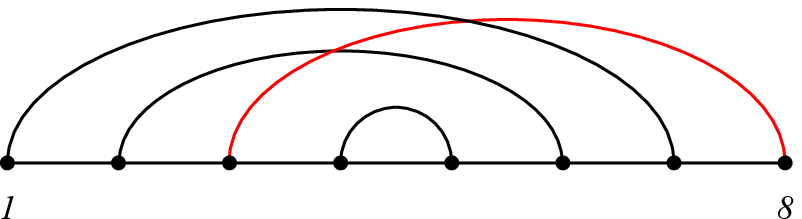}
		\caption{}
	\end{subfigure}
	\begin{subfigure}[b]{0.45\textwidth}
	\centering
    		\includegraphics[scale=0.8]{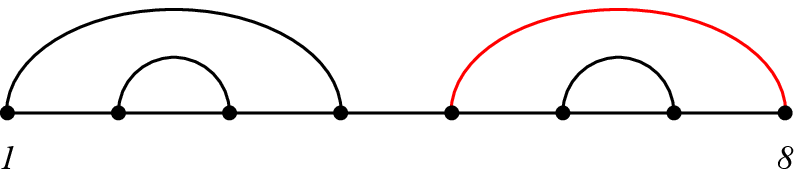}
		\caption{}
	\end{subfigure}
\caption{Counting linear diagrams with counting parallel edges}
\label{fig:simple}
\end{figure}

Let $a_{n,k,l}$ be the number of linear diagrams with $n+1$ chords in total, $k$ loops, and $l$ pairs of parallel chords. We now show that these numbers satisfy the following recurrence relation:
\begin{equation}
\label{eq:a_n_k_l}
a_{n+1,k,l}=a_{n,k-1,l}+(2n+1-k-2l)a_{n,k,l}+(k+1)a_{n,k+1,l}+2(l+1)a_{n,k,l+1}+a_{n,k,l-1}.
\end{equation} 
Take some simple linear diagram with $n+2$ chords and remove a chord ending in the point $2n+4$. Five cases could be possible. The first summand corresponds to the case of removing a single loop $(2n+3,2n+4)$ (figure \ref{fig:simple}(a)). In the second case the numbers of loops and parallel chords don't change (figure \ref{fig:simple}(b), summand $(2n+1-k-2l)a_{n,k,l}$). The third and the fourth cases stand for two situations when a new loop (figure \ref{fig:simple}(c)) or a new pair of parallel chords (figure \ref{fig:simple}(d)), correspondingly, are introduced. Finally, in the fifth case the chord being removed was parallel to some other chord (figure \ref{fig:simple}(e)). 

The recurrence (\ref{eq:a_n_k_l}) has to be solved with the boundary conditions $a_{n,k,l}=0$ when $l>n$, $k>n+1$, $k<0$, or $l<0$, and with the initial conditions $a_{0,k,l}=0$ for all $k$ and $l$, except for $a_{0,1,0}=1$, that corresponds to a diagram with a single chord that is a loop. 

The recurrence (\ref{eq:a_n_k_l}) could be rewritten as  
$$
Q_{n+1,k}(x)=Q_{n,k-1}(x)+(2n+1-k)Q_{n,k}(x)-2x\od{Q_{n,k}(x)}{x}+(k+1)Q_{n,k+1}(x)+xQ_{n,k}(x)+2\od{Q_{n,k}(x)}{x}
$$
for the polynomials 
$$
Q_{n,k}(x)=a_{n,k,0}+a_{n,k,1}x+a_{n,k,2}x^2+\ldots+a_{n,k,n}x^n.
$$
Introducing the polynomials $P_n(z,x)=\sum_{k=0}^nQ_{n,k}(x)z^k$ we get the following equality: 
$$
P_{n+1}(z,x)=(z+x+2n-1)P_n(z,x)-(z-1)\pd{P_n}{z}-2(x-1)\pd{P_n}{x}.
$$
Finally, multiplying this equalitiy by $t^{n}/n!$ and summing over $n$ from $0$ to $+\infty$ yields the following partial differential equation for the generating function $w(t,z,x)=\sum_{n=0}^{+\infty}P_n(z,x)\dfrac{t^{n}}{n!}$:
$$
\pd{w}{t}=(z+x+1)w(z,t)+2t\pd{w}{t}-(z-1)\pd{w}{z}-2(x-1)\pd{w}{x},\qquad\qquad w(0,z,x)=z.
$$
It could be shown that it has the following analytic solution:
$$
w(t,z,x)=\dfrac{(z-1)\sqrt{1-2t}+1}{(1-2t)^{3/2}}\exp\left((z-1)(1-\sqrt{1-2t})-t(1-x)\right).
$$
This function classifies linear diagrams by the number of chords, loops and parallel chords. In particular, setting $z=1$ results in the generating function 
$$
\tilde{w}(t,x)=w(t,1,x)=\dfrac{1}{(1-2t)^{3/2}}\exp\left(t(x-1)\right)=
\sum\limits_{n=0}^{+\infty}\dfrac{t^n}{n!}\sum\limits_{l=0}^n \hat{a}_{n,l}\,x^l,
$$
that classifies linear diagrams by the numbers of chords and pairs of parallel chords only. Its coefficients $\hat{a}_{n,l}$ satisfy the following recurrence relation:
$$
\hat{a}_{n+1,l}=\hat{a}_{n,l-1}+(2n+1-2l)\hat{a}_{n,l}+2(l+1)\hat{a}_{n,l+1}.
$$
Moreover, substituting $z=x=0$ into $w(t,z,x)$, we obtain the generating function $W(t)$ for the numbers $a_{n,0,0}\equiv \bar{a}_n$ enumerating simple linear diagrams with $n+1$ chords:
$$
W(t)=w(t,0,0)=\dfrac{1-\sqrt{1-2t}}{(1-2t)^{3/2}}\exp\left(-1-t+\sqrt{1-2t}\right)=
1\cdot t+3\cdot\dfrac{t^2}{2!}+24\cdot\dfrac{t^3}{3!}+211\cdot\dfrac{t^4}{4!}+
2325\cdot\dfrac{t^5}{5!}+\ldots
$$
From this generating function one could get the following recurrence for the numbers $\bar{a}_n$:
$$
\bar{a}_n=(2n-1)\cdot \bar{a}_{n-1}+(4n-3)\cdot \bar{a}_{n-2}+(2n-4)\cdot \bar{a}_{n-3},\qquad n\geq 2;\qquad \text{$\bar{a}_n=0$ for $n\leq 0$;}\quad \bar{a}_1=1.
$$

Now we turn to enumerating simple chord diagrams. As a result of gluing a simple linear diagram with $n$ chords we get a simple chord diagram with $n$ chords if the original linear diagram neither had a chord $(1,2n)$ nor a pair of chords $(1,i)$, $(i+1,2n)$ for some $i$. Let $\bar{q}_n$ be the number of diagrams with $n$ chords and no chord $(1,2n)$. We claim that 
$$
\bar{q}_n=\bar{a}_{n-1}-\bar{q}_{n-1}.
$$
Consider a simple diagram with a chord $(1,2n)$. Removing this chord, we obtain a simple diagram with $n-1$ chords and no chord $(1,2n-2)$. Indeed, the existence of such a chord would mean that the original diagram had two parallel chords, but that is impossible by definition of a simple diagram. Thus we conclude that the number of simple diagrams with the chord $(1,2n)$ is equal to $\bar{q}_{n-1}$. 

Let us now prove that the number of simple chord diagrams $\bar{b}_n$ with $n$ chords satisfies the recurrence 
$$
\bar{b}_n=\bar{q}_n-\bar{b}_{n-1}.
$$
Indeed, gluing any linear diagram without a chord $(1,2n)$ into a chord diagram we get a pair of parallel chords if and only if the original diagram had chords $(1,i)$ and $(i+1,2n)$ for some $i$. Removing the chord $(i+1,2n)$ we get a linear diagram without such pairs and without a chord $(1,2n-2)$. The number of such linear diagrams is equal to $\bar{b}_{n-1}$. 

It terms of exponential generating functions $W(t)$, $V(t)$, and $U(t)$ for the numbers $\bar{a}_n$, $\bar{q}_n$ and $\bar{b}_n$ we get the system 
$$
\begin{aligned}
V'(t)&=W(t)-V(t),\\
U'(t)&=V'(t)-U(t).
\end{aligned}
$$
Solving this system yields the following exponential generating function $U(t)$ for labelled simple chord diagrams:
$$
U(t)=\dfrac{e^{-1-t+\sqrt{1-2t}}}{\sqrt{1-2t}}(1+\sqrt{1-2t})-(2-t)e^{-t}=1\cdot \dfrac{t^2}{2!}+1\cdot \dfrac{t^3}{3!}+
21\cdot \dfrac{t^4}{4!}+168\cdot \dfrac{t^5}{5!}+1968\cdot \dfrac{t^6}{6!}+\ldots
$$

Let us turn to unlabelled simple diagrams. The technique of enumerating them doesn't differ much from the one in the previous section. However, for these diagrams one has to consider more different cases to obtain the recurrence relations. In order not to overload the article, we will still describe all of them, but maybe omit some details about proving the exact forms for the coefficients and for the initial conditions. It should be straightforward to check them.

\begin{figure}[ht]
\centering
	\begin{subfigure}[b]{0.24\textwidth}
	\centering
    		\includegraphics[scale=0.8]{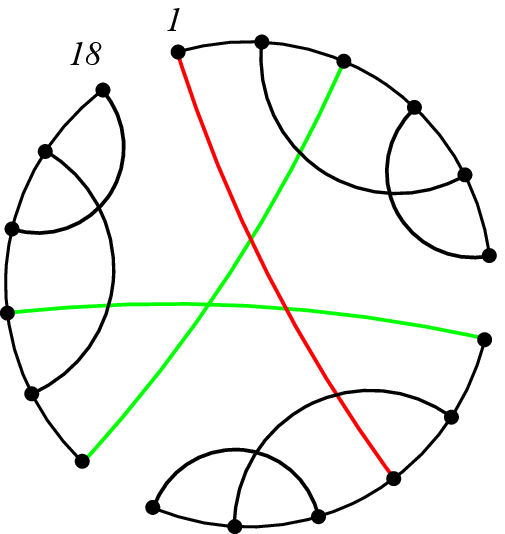}
		\caption{}
	\end{subfigure}
	\begin{subfigure}[b]{0.24\textwidth}
	\centering
    		\includegraphics[scale=0.8]{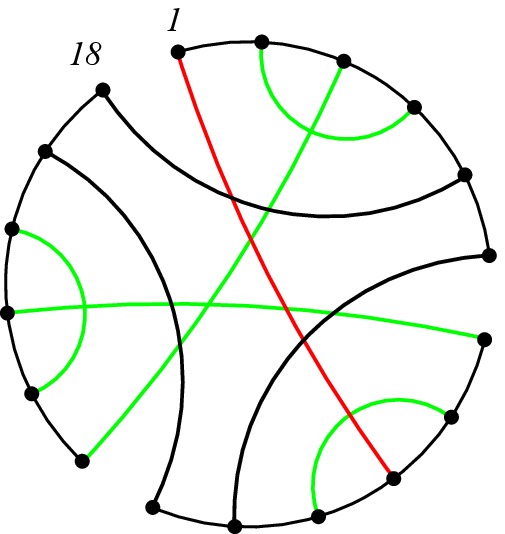}
		\caption{}
	\end{subfigure}	
	\begin{subfigure}[b]{0.24\textwidth}
	\centering
    		\includegraphics[scale=0.8]{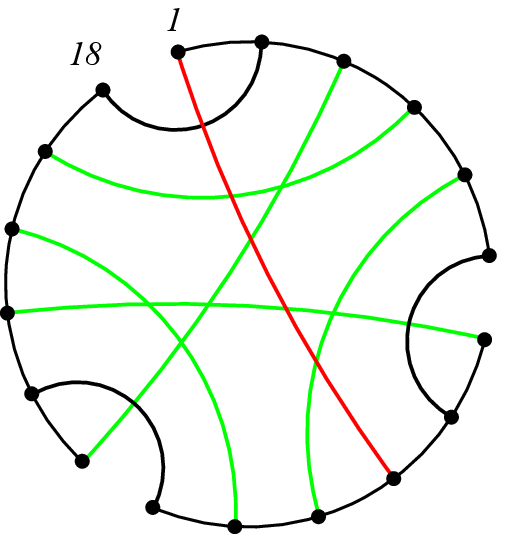}
		\caption{}
	\end{subfigure}	
	\begin{subfigure}[b]{0.24\textwidth}
	\centering
    		\includegraphics[scale=0.8]{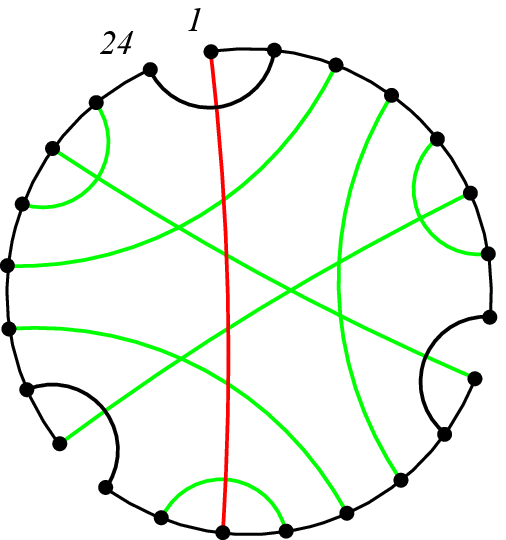}
		\caption{}
	\end{subfigure}	
	\caption{Counting simple $d$-symmetric linear diagrams, $d$ odd}
\label{fig:d_odd_symmetric_linear}
\end{figure}

Let $\bar{a}_{m}^{(d)}$ be the numbers of $d$-symmetric linear simple diagrams with $m\cdot d$ points. In the case of odd $d$ these numbers satisfy the following recurrence:
$$
\bar{a}_{m}^{(d)}=[(m-1)\cdot d-2]\cdot \bar{a}_{m-2}^{(d)}+(2m-7)\cdot d\cdot \bar{a}_{m-4}^{(d)}+(m-6)\cdot d\cdot \bar{a}_{m-6}^{(d)} \quad 
\text{for $m>2$;}
$$
$$
\bar{a}_{m}^{(d)}=0 \quad\text{for $m<0$ or $m=1$;}\quad \bar{a}_{0}^{(d)}=1;\quad \bar{a}_{2}^{(d)}=d-1.
$$
To prove this formula we consider four cases that could be possible for a chord $\{1,i\}$:
\begin{itemize}
\item After removing this chord the diagram is still simple (figure \ref{fig:d_odd_symmetric_linear}(a)).
\item Removing this chord creates a loop, but after removing this loop the diagram becomes simple (figure \ref{fig:d_odd_symmetric_linear}(b)). 
\item Removing a chord $\{1,i\}$ creates a pair of parallel chords. If that happens, we remove one (arbitrary) chord from this pair too (figure \ref{fig:d_odd_symmetric_linear}(c)). 
\item We remove an chord $\{1,i\}$, a loop appears, but removing this loop results in a pair of parallel chords (figure \ref{fig:d_odd_symmetric_linear}(d)). In this case we remove one chord from that pair too.
\end{itemize}

\begin{figure}[ht]
\centering
	\begin{subfigure}[b]{0.3\textwidth}
	\centering
    		\includegraphics[scale=0.8]{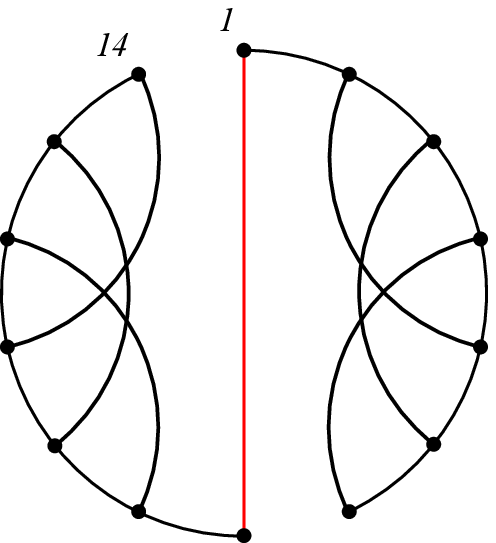}
		\caption{}
	\end{subfigure}	
	\begin{subfigure}[b]{0.3\textwidth}
	\centering
    		\includegraphics[scale=0.8]{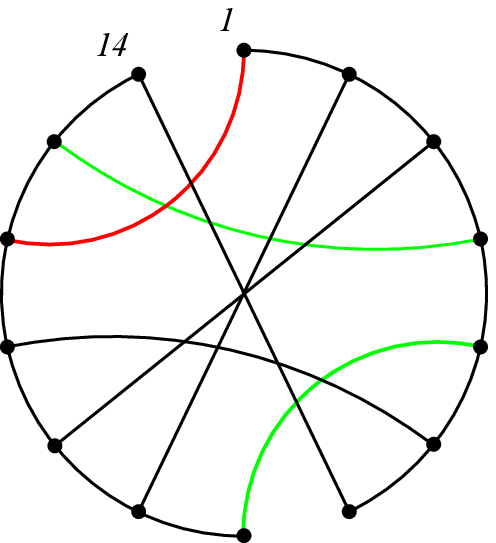}
		\caption{}
	\end{subfigure}
	\begin{subfigure}[b]{0.3\textwidth}
	\centering
    		\includegraphics[scale=0.8]{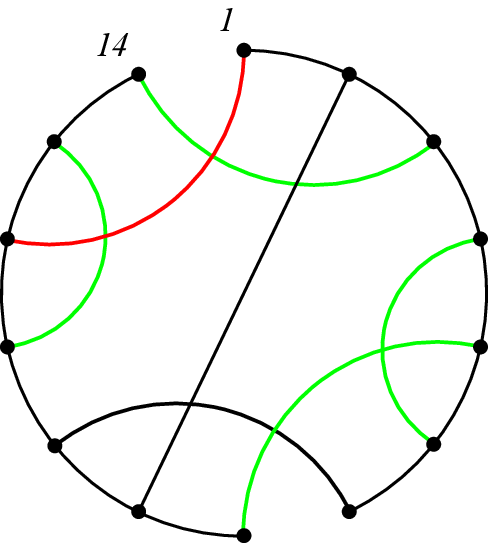}
		\caption{}
	\end{subfigure}
\caption{Counting simple $d$-symmetric linear diagrams, $d$ even}
\label{fig:d_symmetric_even}
\end{figure}

In the case of even $d$ the situation is a bit more complicated. Namely, we have to introduce a new parameter $k$. Let $\bar{a}_{m,k}^{(d)}$ be the number of $d$-symmetric linear diagrams with $m\cdot d$ chords, $k\cdot d/2$ of which connect the opposite points of a diagram. These number equal 
$$
\bar{a}_{m,k}^{(d)}=0\quad \text{for $m<0$ or $k<0$ or $k>m$,}\qquad \bar{a}_{0,0}^{(d)}=1,
$$
and satisfy the recurrence
$$
\bar{a}_{m,k}^{(d)}=\bar{a}_{m-1,k-1}^{(d)}
+[(m-1)\cdot d -3+\delta_{m,2}]\cdot \bar{a}_{m-2,k}^{(d)}
+(k+1)\cdot d\cdot \bar{a}_{m-3,k+1}^{(d)}+
$$
$$
+(2m+k-7)\cdot \bar{a}_{m-4,k}^{(d)}
+(k+1)\cdot d\cdot \bar{a}_{m-5,k+1}^{(d)}
+(m+k-6)\cdot d\cdot \bar{a}_{m-6,k}^{(d)}
$$
for other values of $m,k$. Here $\delta_{m,2}$ is the Kronecker delta that equals $1$ if $m=2$ and $0$ otherwise. Three types of diagrams that are new here compared to the odd case are shown on figure \ref{fig:d_symmetric_even}. The summand $\bar{a}_{m-1,k-1}^{(d)}$ enumerates diagrams with a chord $\{1, i\}$ that is a diameter (figure \ref{fig:d_symmetric_even}(a)). The summand $(k+1)\cdot d\cdot \bar{a}_{m-3,k+1}^{(d)}$ corresponds to the case when removing a chord $\{1, i\}$ creates $d/2$ pairs of parallel chords, and we remove one chord from each pair, so that the other one becomes a diameter (figure \ref{fig:d_symmetric_even}(b)). The summand $(k+1)\cdot d\cdot \bar{a}_{m-5,k+1}^{(d)}$ describes a similar case with an intermediate step of removing $d$ loops (figure \ref{fig:d_symmetric_even}(c)). 

\begin{figure}[ht]
\centering
	\begin{subfigure}[b]{0.3\textwidth}
	\centering
    		\includegraphics[scale=0.8]{pics/simple_linear_2_symmetric__loops.eps}
		\caption{}
	\end{subfigure}	
	\begin{subfigure}[b]{0.3\textwidth}
	\centering
    		\includegraphics[scale=0.8]{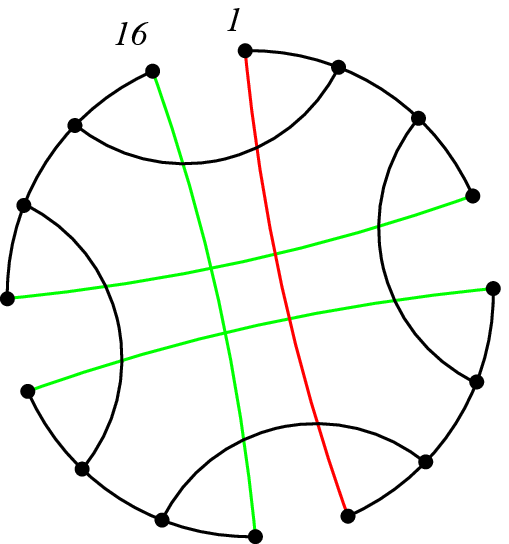}
		\caption{}
	\end{subfigure}
	\begin{subfigure}[b]{0.3\textwidth}
	\centering
    		\includegraphics[scale=0.8]{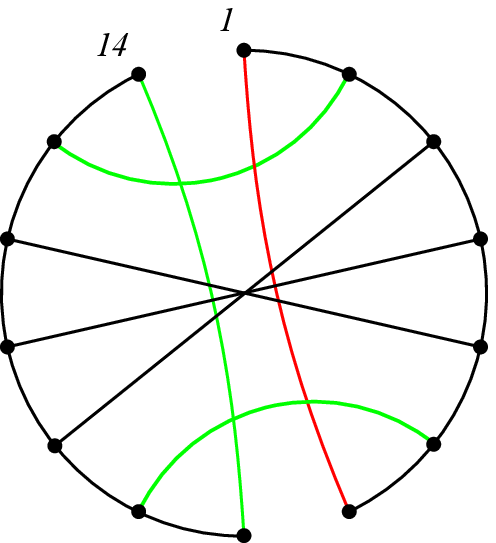}
		\caption{}
	\end{subfigure}
\caption{Special kinds of $d$-symmetric linear diagrams}
\label{fig:2_symmetric}
\end{figure}

For both odd and even cases the numbers $\bar{a}_m^{(d)}$ are obviously equal to $\sum_{k=0}^m\bar{a}_{m,k}^{(d)}$. We turn to enumerating unlabelled simple chord diagrams. First we define the auxiliary sequence $\bar{q}_m^{(d)}$ as follows:
$$
\bar{q}_m^{(d)}=\bar{a}_m^{(d)}-\bar{q}_{m-2}^{(d)}\qquad \text{for $m\geq 2$};\qquad
\bar{q}_0^{(d)}=\bar{q}_1^{(d)}=\begin{cases}
\,0,&\quad \text{$d\leq 2$,}\\
\,1,&\quad \text{$d>2$;}
\end{cases}\qquad \bar{q}_m^{(d)}=0,\,\,m<0.
$$
We claim that $\bar{q}_m^{(d)}$ enumerates $d$-symmetric loopless chord diagrams with $m\cdot d$ points that could be obtained by gluing simple $d$-symmetric linear diagrams. Indeed, a loop can appear after gluing only if the vertices $1$ and $m\cdot d$ were joined by a chord (figure \ref{fig:2_symmetric}(a)). If a diagram has such a chord, removing it and all the chords on its orbit results in a diagram having $(m-2) \cdot d$ points that cannot have chords of the same type anymore, as this would mean that the original diagram had parallel chords. By induction, such diagrams are enumerated by $\bar{q}_{m-2}^{(d)}$. Another class of linear $d$-symmetric diagrams enumerated by $\bar{q}_m^{(d)}$ if $d=2k$ for $k \geq 2$ are simple diagrams with $m \cdot d$ points that don't have a chord $\{1,m\cdot d/2\}$. This can be proved by applying the same argument with chord removal for the diagrams that have this chord (figure \ref{fig:2_symmetric}(b)).

To proceed, we also have to examine the case when a $2$-symmetric diagram does neither have a chord $\{1,m\}$ nor a chord $\{1,2m\}$. The corresponding numbers $\bar{p}_m$ are given by
$$
\bar{p}_m=\bar{q}_m^{(2)}-\bar{p}_{m-2}-\bar{q}_{m-4}^{(2)},\qquad \text{for $m\geq 2$};\qquad \bar{p}_0=0,\quad \bar{p}_1=1.
$$
Indeed, the summand $\bar{q}_m^{(2)}$ enumerates those diagrams that don't have a chord $\{1,2m\}$. A subset of those diagrams that have a chord $\{1,m\}$ is enumerated by $\bar{p}_{m-2}+\bar{q}_{m-4}^{(2)}$: there are $\bar{p}_{m-2}$ diagrams without a chord $\{2,2m-1\}$ and $\bar{q}_{m-4}^{(2)}$ diagrams with such a chord (figure \ref{fig:2_symmetric}(c)).

\begin{figure}[ht]
\centering
	\begin{subfigure}[b]{0.3\textwidth}
	\centering
    		\includegraphics[scale=0.8]{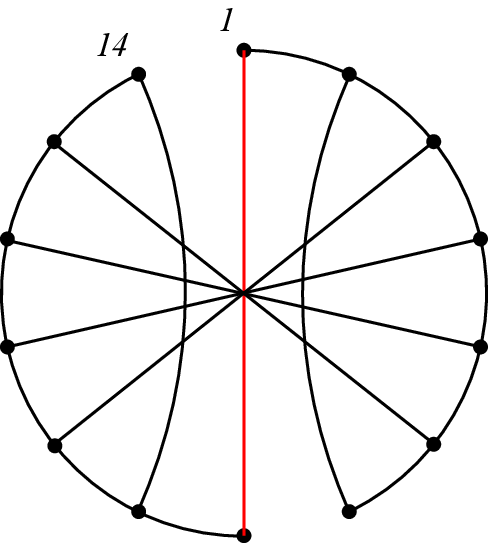}
		\caption{}
	\end{subfigure}
	\begin{subfigure}[b]{0.3\textwidth}
	\centering
    		\includegraphics[scale=0.8]{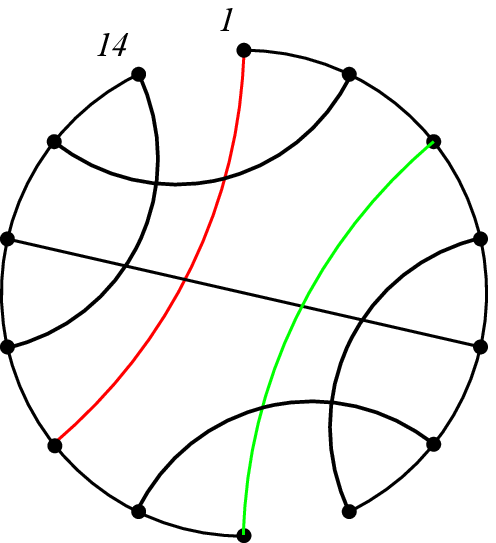}
		\caption{}
	\end{subfigure}	
	\caption{Counting $2$-symmetric chord diagrams}
\label{fig:2_symmetric_parallel}
\end{figure}

Now we are ready to give the recurrence for the number of $2$-symmetric simple chord diagrams with $2m$ points:
$$
\bar{f}(2m,2)=\bar{p}_m-\bar{f}(2m-4,2)-\bar{q}_{m-1}^{(2)}+\bar{p}_{m-1}\qquad \text{for $m\geq 2$};\qquad \bar{f}(0,2)=1,\quad \bar{f}(2,2)=0;
$$
In this equation we start with the number $\bar{p}_m$ of diagrams that are both loopless and free of chords $\{1,m\}$ and $\{m+1,2m\}$, and subtract the number of those diagrams that have other kinds of parallel chords. Observe that the diagrams enumerated by $\bar{p}_m$ are obtained from simple linear diagrams, so if they have parallel chords, one pair of them has the form $(\{1,i\}, \{2m,i+1\})$ and the other pair is just obtained from the first one by a rotation by $180^{\circ}$. Two cases could be possible: either the chord $\{1,i\}$ is a diameter, that is $i=m+1$, or not. In the first case (figure \ref{fig:2_symmetric_parallel}(a)) removing the diameter results in a simple $2$-symmetric diagram with $2m-2$ points and a pair of chords $(\{1,m-1\}, \{m,2m-2\}$). These diagrams are enumerated by $\bar{q}_{m-1}^{(2)}-\bar{p}_{m-1}$. In the second case (figure \ref{fig:2_symmetric_parallel}(b)) by removing the chord $\{1,i\}$ together with the one symmetric to it we get a linear diagram that is ready to be glued into a simple chord diagram, and the corresponding summand is $\bar{f}(2m-4,2)$.

For $d$-symmetric simple chord diagrams where $d>2$ we have the following enumeration formulas:
\begin{equation}
\label{eq:d_symmetric_chord_odd}
\bar{f}((2i+1)\cdot m,2i+1)=\bar{q}_m^{(2i+1)}-\bar{f}((2i+1)\cdot (m-2),2i+1)\quad \text{for $m\geq 1$};\quad \bar{f}(0,2i+1)=0;\quad i\geq 0;
\end{equation}
\begin{equation}
\label{eq:d_symmetric_chord_even}
\bar{f}(2i\cdot m,2i)=\bar{q}_m^{(2i)}-\bar{f}(2i\cdot (m-2),2i)-\bar{q}_{m-2}^{(2i)}-\bar{q}_{m-3}^{(2i)}\qquad \text{for $m\geq 1$};\qquad \bar{f}(0,2i)=0;\qquad i>1.
\end{equation}
\begin{figure}[ht]
\centering
	\begin{subfigure}[b]{0.3\textwidth}
	\centering
    		\includegraphics[scale=0.8]{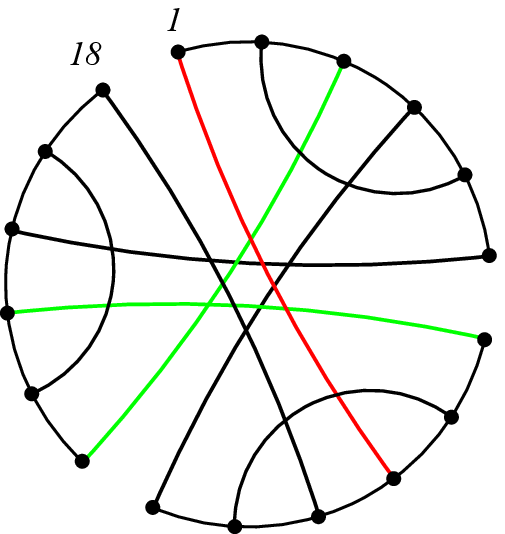}
		\caption{}
	\end{subfigure}
	\begin{subfigure}[b]{0.3\textwidth}
	\centering
    		\includegraphics[scale=0.8]{pics/simple_chord_d_symmetric__even_1.eps}
		\caption{}
	\end{subfigure}	
	\caption{Counting $d$-symmetric chord diagrams, $d>2$}
\label{fig:d_symmetric_chord}
\end{figure}
It's easy to prove (\ref{eq:d_symmetric_chord_odd}): in addition to the possibility of loops which was already excluded in $\bar{q}_m^{(2i+1)}$, we also exclude those diagrams that contain $d$ pairs of parallel chords after gluing (figure \ref{fig:d_symmetric_chord}(a)). Clearly, those are enumerated by $\bar{f}(2i\cdot (m-2),2i)$: removing a chord $\{1,j\}$ establishes a bijection with the corresponding type of diagrams. In $(\ref{eq:d_symmetric_chord_even})$ two additional summands stand for two special cases when removing a chord $\{1,j\}$ does not lead to an arbitrary $2i$-symmetric simple chord diagram with less chords. The first case was already shown on figure \ref{fig:2_symmetric}(b), and corresponds to the summand $\bar{q}_{m-2}^{(2i)}$. The second one is the case of a chord $\{1,i\cdot d\}$ being parallel to two different chords (figure \ref{fig:d_symmetric_chord}(b)). Removing this chord transforms this case into the first one, so the number of corresponding diagrams is $\bar{q}_{m-3}^{(2i)}$.

All the numbers in the right hand side of the Burnside lemma (\ref{eq:Burnside_lemma_Cn}) are now known, and that gives us the following sequence enumerating simple chord diagrams not isomorphic under rotations: $0,1,1,4,21,176,1893,\ldots$ (see Table \ref{table:simple}).

Finally, we will enumerate unlabelled simple chord diagrams for the case when our notion of isomorphism includes reflection symmetries as well. As before, we begin with enumerating some objects analogous to linear diagrams, as it is more convenient to establish a recurrence for them. 

\begin{figure}[ht]
\centering
	\begin{subfigure}[b]{0.24\textwidth}
	\centering
    		\includegraphics[scale=0.8]{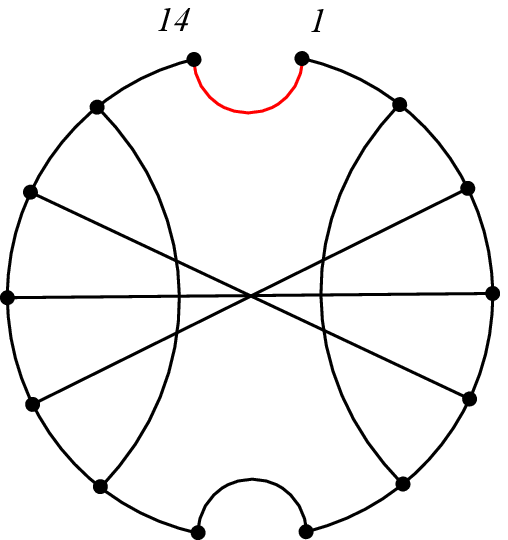}
		\caption{}
	\end{subfigure}
	\begin{subfigure}[b]{0.24\textwidth}
	\centering
    		\includegraphics[scale=0.8]{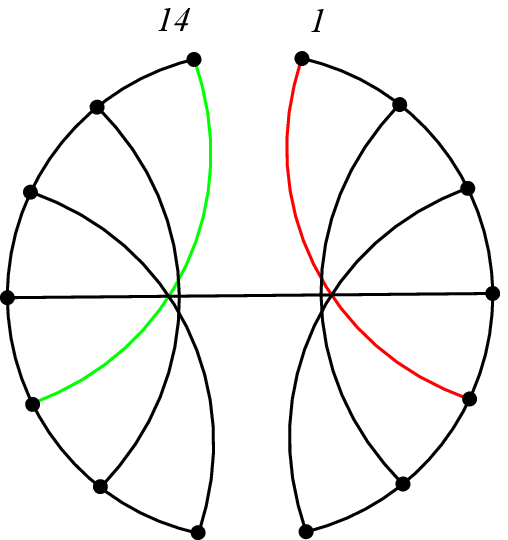}
		\caption{}
	\end{subfigure}
	\begin{subfigure}[b]{0.24\textwidth}
	\centering
    		\includegraphics[scale=0.8]{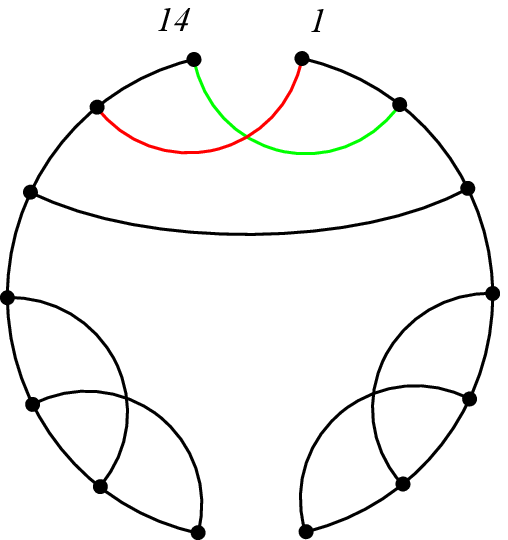}
		\caption{}
	\end{subfigure}
	\begin{subfigure}[b]{0.24\textwidth}
	\centering
    		\includegraphics[scale=0.8]{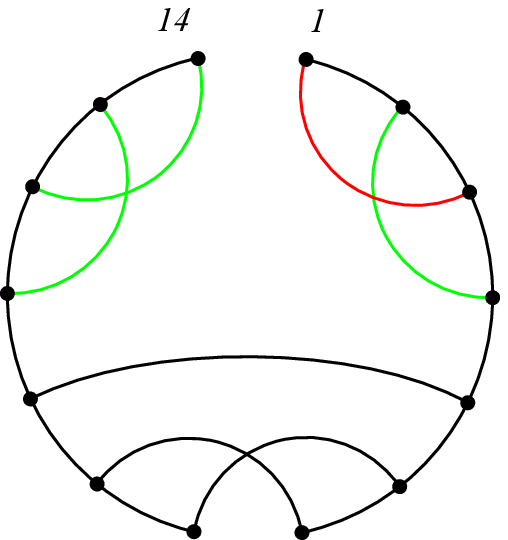}
		\caption{}
	\end{subfigure}
	\caption{Simple linear diagrams with reflectional symmetry}
\label{fig:refl_symmetric_chord}
\end{figure}

Let $\bar{r}_{n,k}$ be the number of simple diagrams with the following properties: they have $2n$ points, points $1$ and $2n$ as well as $n$ and $n+1$ are not thought to be neighboring, reflection across the line passing through the center of a diagram and the center of the segment connecting points $1$ and $2n$ transforms a diagram into itself, and the number of chords that are transformed into itself under such a reflection is equal to $k$. For example the diagram on figure \ref{fig:refl_symmetric_chord}(a) is the one counted by $\bar{r}_{7,3}$. For the number of such diagrams with an additional property that they have a chord $\{1,2n\}$ we introduce a separate sequence $\bar{s}_{n,k}$. We claim that

$$
\bar{s}_{n,k} = \bar{r}_{n-1,k-1} - \bar{s}_{n-1,k-1} \text{ for } n > 0; \qquad \bar{s}_{0,k}=0.
$$
Indeed, any corresponding diagram consists of a chord $\{1,2n\}$ and the remaining part which is a diagram with $2n-2$ points that does not have the same type of chord, thus counted by $\bar{r}_{n-1,k-1} - \bar{s}_{n-1,k-1}$. The recurrence for $\bar{r}_{n,k}$ is much more complicated:
$$
\bar{r}_{n,k} = 0 \quad \text{ if } \quad n < 0 \text{ or } k < 0 \text{ or } k > n.  \qquad \bar{r}_{0,0} = \bar{r}_{2,0} = 1.
$$
$$
\bar{r}_{n,k} = \bar{s}_{n,k} + 2(n-2)\bar{r}_{n-2,k} + \bar{s}_{n-2,k} + 2(2n-k-7)\bar{r}_{n-4,k} + 2(k-1)\bar{r}_{n-3,k-1}+ $$
$$
+ 2(k-1)\bar{r}_{n-5,k-1} + 2(n-k-6)\bar{r}_{n-6,k}\quad  \text { otherwise}.
$$

\begin{figure}[ht]
\centering
	\begin{subfigure}[b]{0.24\textwidth}
	\centering
    		\includegraphics[scale=0.8]{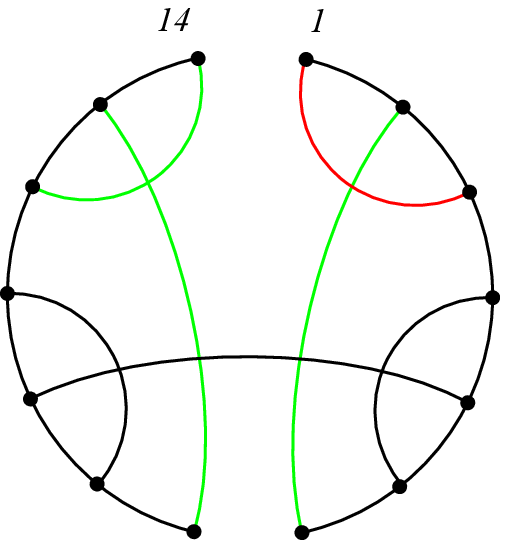}
		\caption{}
	\end{subfigure}
	\begin{subfigure}[b]{0.24\textwidth}
	\centering
    		\includegraphics[scale=0.8]{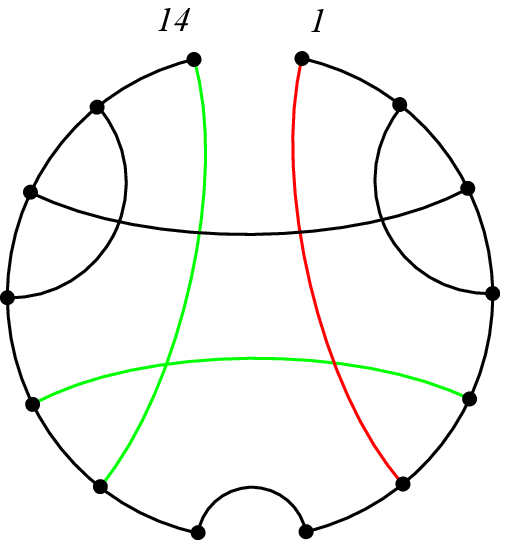}
		\caption{}
	\end{subfigure}
	\begin{subfigure}[b]{0.24\textwidth}
	\centering
    		\includegraphics[scale=0.8]{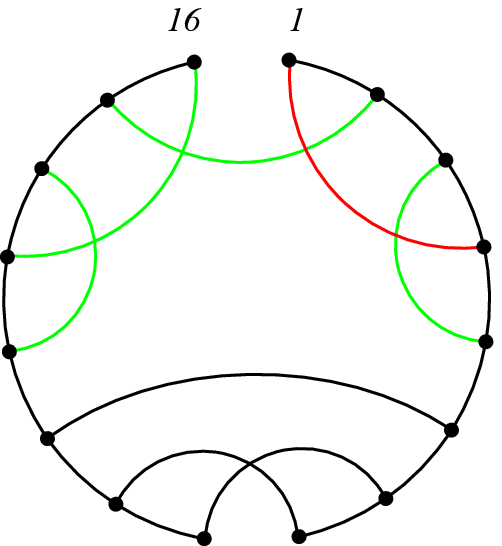}
		\caption{}
	\end{subfigure}
	\begin{subfigure}[b]{0.24\textwidth}
	\centering
    		\includegraphics[scale=0.8]{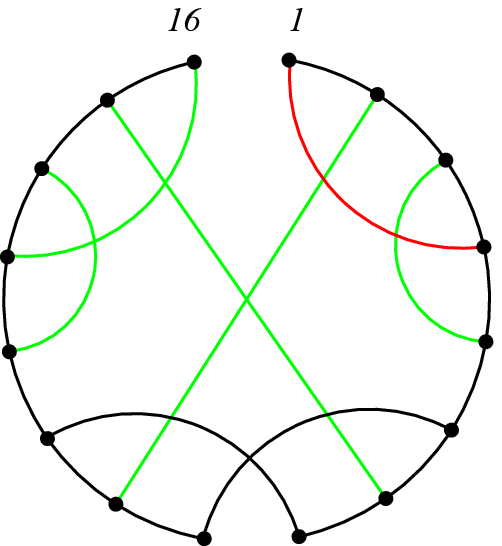}
		\caption{}
	\end{subfigure}
	\caption{Counting simple diagrams with reflectional symmetry}
\label{fig:refl_symmetric_chord_2}
\end{figure}

To prove it we again classify all the corresponding diagrams according to the properties of a chord $\{1,i\}$. 7 possible cases are:
\begin{itemize}
\item $i = 2n$. This case is counted by $\bar{s}_{n,k}$. In all the cases below it is assumed that $i \neq 2n$.
\item Removing a chord $\{1,i\}$ together with the one symmetric to it yields a simple diagram (figure \ref{fig:refl_symmetric_chord}(b)). This case is enumerated by $2(n-2)\bar{r}_{n-2,k} + \bar{s}_{n-2,k}$, because to reconstruct an initial diagram one could take a diagram with $2n-4$ points and find a place for the new chord. $2n-4$ options are always available, and one more option becomes available if the chord $\{1,2n-4\}$ is present, because it cannot be parallel to the newly added chord (result is shown on figure \ref{fig:refl_symmetric_chord}(c)).
\item Removing a chord $\{1,i\}$ together with the one symmetric to it yields a diagram with a pair of loops; removing these loops yields a simple diagram. (figure \ref{fig:refl_symmetric_chord}(d), summand $2(n-3)\bar{r}_{n-4,k}$).

\item Removing a chord $\{1,i\}$ together with the one symmetric to it yields a diagram with two pairs of parallel chords; removing one chord from each pair yields a simple diagram (figure \ref{fig:refl_symmetric_chord_2}(a), summand $2(n-k-4)\bar{r}_{n-4,k}$). 

\item Removing a chord $\{1,i\}$ together with the one symmetric to it yields a diagram with one pair of parallel chords; removing one chord from it yields a simple diagram. (figure \ref{fig:refl_symmetric_chord_2}(b), summand $2(k-1)\bar{r}_{n-3,k-1}$).

\item Removing a chord $\{1,i\}$ together with the one symmetric to it yields a diagram with two pairs of loops; removing them creates one pair of parallel chords; removing a chord from this pair yields a simple diagram. (figure \ref{fig:refl_symmetric_chord_2}(c), summand $2(k-1)\bar{r}_{n-5,k-1}$).

\item Removing a chord $\{1,i\}$ together with the one symmetric to it yields a diagram with two pairs of loops; removing them creates two pairs of parallel chords; removing a chord from each pair yields a simple diagram. (figure \ref{fig:refl_symmetric_chord_2}(c), summand $2(n-k-6)\bar{r}_{n-6,k}$).

\end{itemize}

\begin{figure}[ht]
\centering
	\begin{subfigure}[b]{0.3\textwidth}
	\centering
    		\includegraphics[scale=0.8]{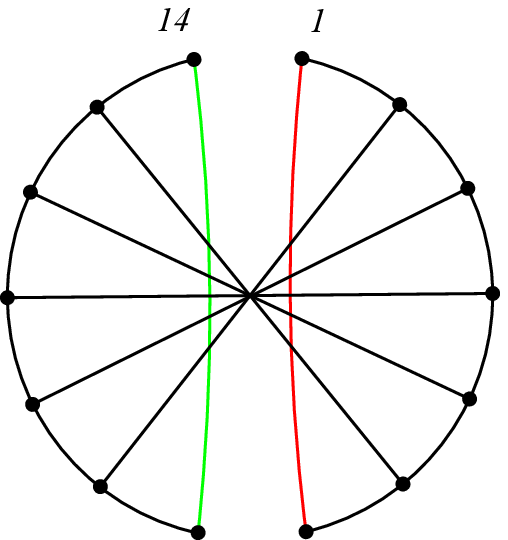}
		\caption{}
	\end{subfigure}
	\begin{subfigure}[b]{0.3\textwidth}
	\centering
    		\includegraphics[scale=0.8]{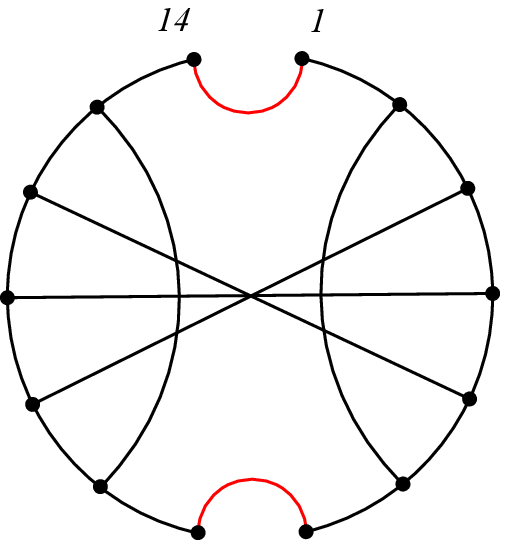}
		\caption{}
	\end{subfigure}	
	\caption{Counting chord diagrams with reflectional symmetry}
\label{fig:refl_symmetric_chord_3}
\end{figure}

Denoting by $\bar{r}_n$ the sum of $\bar{r}_{n,k}$ over $k$ and by $\bar{s}_n$ the corresponding sum of $\bar{s}_{n,k}$, we are ready to enumerate chord diagrams that are transformed into themselves by reflections. Those can be of two types. The axis of symmetry can either pass through two opposite points, or through the center of the segment that connects two neighboring points, like in linear diagrams we studied above. Let the number of diagrams of the first type with $n$ chords be $\bar{K(n)}$, and of the second type $\bar{H}(n)$. We claim that
$$
\bar{K}(0) = \bar{K}(1) = 0; \qquad \bar{K}(n) = \bar{r}_{n-1} - \bar{K}(n-2) \text{ for } n > 1.
$$
Indeed, any diagram enumerated by $K_n$ has a chord on its axis of symmetry. Removing this chord yields a diagram that can be almost arbitrary, except that it can't have chords $\{1,n-1\}$ and $\{n,2n-2\}$ like the one on figure \ref{fig:refl_symmetric_chord_3} (a). But this forbidden type of diagrams is enumerated by $\bar{K}(n-2)$, as these two chords can be viewed as a single chord that lies on the axis of symmetry. For $\bar{H}(n)$ the formula is
$$
\bar{H}(0) = 0; \qquad \bar{H}(n) = \bar{L}(n) - \bar{K}(n-1)  \text{ for } n > 0, \text{ where } 
$$
$$
\bar{L}(0) = \bar{L}(1) = 0; \qquad \bar{L}(n)=\bar{r}_n-2\bar{s}_n+\bar{L}(n-2) \text{ for } n > 1.
$$

Here $\bar{K}(n-1)$ is the number of chord diagrams that would have parallel chords after gluing (figure \ref{fig:refl_symmetric_chord_3}(a)). It is subtracted from $\bar{L}(n)$, which is the number of diagrams that would not have any loops. The given recurrence for $\bar{L}(n)$ is easily proved: take all $\bar{r}_n$ diagrams, subtract those with a chord $\{1,2n\}$, then those with a chord $\{n, n+1\}$, and finally add the number of diagrams that have both of them. The last number is $\bar{L}(n-2)$, as deleting these chords (figure \ref{fig:refl_symmetric_chord_3}(b)) yields a diagram that doesn't have this type of chords anymore. 

The following formula enumerates non-isomorphic simple diagrams under the action of a dihedral group. It is just the Burnside lemma rewritten in terms of the sequences that we derived above:
\begin{equation}
\label{eq:Burnside_lemma_Dn_simple}
\bar{c}_n=\dfrac{1}{4n}\left[\sum\limits_{d\,|\, 2n}\phi(d)\,\bar{f}(2n,d)+n\cdot \bar{K}(n)+n\cdot \bar{H}(n)\right].
\end{equation} 

The sequence obtained from it starts with the numbers $0,1,1,4,18,116,1060$ (see Table \ref{table:simple}).

\section*{Conclusion}
In the first part of this paper labelled and unlabelled loopless chord diagrams were enumerated.  Consequently, we obtained the expressions for the numbers of non-isomorphic hamiltonian cycles in an unlabelled $n$-dimensional octahedron. In the second part we solved a more technically complex problem of enumerating simple diagrams. We obtained a generating function for labelled diagrams classified by the numbers of loops and parallel chords. As a special case, we got a generating function for simple chord diagrams. For unlabelled simple chord diagrams we gave a solution in a form of several recurrence relations. The following tables list the numbers of different kinds of diagrams described above.

\begin{table}[h!]
\scriptsize
\centering
\begin{tabular}{c|cccc}
\midrule
$ n $  &\phantom{0}Linear\phantom{0}&\phantom{0}Chord labelled\phantom{0}&\phantom{0}Under rotations\phantom{0}&\phantom{0}Under all symmetries\phantom{0}\\ 
\midrule
1 & 0 & 0 & 0 & 0 \\
2 & 1 & 1 & 1 & 1 \\
3 & 5 & 4 & 2 & 2 \\
4 & 36 & 31 & 7 & 7 \\
5 & 329 & 293 & 36 & 29 \\
6 & 3655 & 3326 & 300 & 196 \\
7 & 47844 & 44189 & 3218 & 1788 \\
8 & 721315 & 673471 & 42335 & 21994 \\
9 & 12310199 & 11588884 & 644808 & 326115 \\
10 & 234615096 & 222304897 & 11119515 & 5578431 \\
11 & 4939227215 & 4704612119 & 213865382 & 107026037 \\
12 & 113836841041 & 108897613826 & 4537496680 & 2269254616 \\
13 & 2850860253240 & 2737023412199 & 105270612952 & 52638064494 \\
14 & 77087063678521 & 74236203425281 & 2651295555949 & 1325663757897 \\
15 & 2238375706930349 & 2161288643251828 & 72042968876506 & 36021577975918 \\
16 & 69466733978519340 & 67228358271588991 & 2100886276796969 & 1050443713185782 \\
17 & 2294640596998068569 & 2225173863019549229 & 65446290562491916 & 32723148860301935 \\
18 & 80381887628910919255 & 78087247031912850686 & 2169090198219290966 & 1084545122297249077 \\
19 & 2976424482866702081004 & 2896042595237791161749 & 76211647261082309466 & 38105823782987999742 \\
20 & 116160936719430292078411 & 113184512236563589997407 & 2829612806029873399561 & 1414806404051118314077 \\
\midrule\end{tabular}
\caption{Loopless diagrams by number of chords}
\label{table:loopless}
\end{table}

\begin{table}[h!]
\scriptsize
\centering
\begin{tabular}{c|cccc}
\midrule
$ n $  &\phantom{0}Linear\phantom{0}&\phantom{0}Chord labelled\phantom{0}&\phantom{0}Under rotations\phantom{0}&\phantom{0}Under all symmetries\phantom{0}\\ 
\midrule
1 & 0 & 0 & 0 & 0 \\
2 & 1 & 1 & 1 & 1 \\
3 & 3 & 1 & 1 & 1 \\
4 & 24 & 21 & 4 & 4 \\
5 & 211 & 168 & 21 & 18 \\
6 & 2325 & 1968 & 176 & 116 \\
7 & 30198 & 26094 & 1893 & 1060 \\
8 & 452809 & 398653 & 25030 & 13019 \\
9 & 7695777 & 6872377 & 382272 & 193425 \\
10 & 146193678 & 132050271 & 6604535 & 3313522 \\
11 & 3069668575 & 2798695656 & 127222636 & 63667788 \\
12 & 70595504859 & 64866063276 & 2702798537 & 1351700744 \\
13 & 1764755571192 & 1632224748984 & 62778105236 & 31390695708 \\
14 & 47645601726541 & 44316286165297 & 1582725739329 & 791372281393 \\
15 & 1381657584006399 & 1291392786926821 & 43046433007765 & 21523271532811 \\
16 & 42829752879449400 & 40202651019430461 & 1256332883208474 & 628166776833181 \\
17 & 1413337528735664887 & 1331640833909877144 & 39165907107963273 & 19582955637428422 \\
18 & 49465522112961344241 & 46762037794122159492 & 1298945495674093932 & 649472761243051940 \\
19 & 1830184115528550306438 & 1735328399106396110310 & 45666536827274985585 & 22833268501579122332 \\
20 & 71375848864779552073957 & 67858430028772637693845 & 1696460750775267473762 & 848230375982060558217 \\
\midrule\end{tabular}
\caption{Simple diagrams by number of chords}
\label{table:simple}
\end{table}

\newpage

\bibliography{bib}

\end{document}